\documentclass[reqno,11pt]{amsart}
\usepackage{amsmath,amssymb,mathrsfs,amsthm,amsfonts}
\usepackage[usenames,dvipsnames]{xcolor}
\usepackage{comment}
\usepackage{hyperref}
\usepackage{graphicx}
\usepackage{enumerate}
\usepackage{stmaryrd}
\usepackage[all]{xy}
\usepackage{bbm}
\usepackage{hyperref}
\usepackage{caption}
\usepackage{subcaption}
\hypersetup{
     colorlinks=true,
     linkcolor=blue,
     filecolor=black,
     citecolor = ForestGreen,      
     urlcolor=blue,
     }
\usepackage{fullpage}

\usepackage{tikz}
\usetikzlibrary{decorations.markings}
\usetikzlibrary{arrows}
\usetikzlibrary{calc}

\newtheorem{theorem}{Theorem}[section]
\newtheorem{lemma}[theorem]{Lemma}
\newtheorem{claim}[theorem]{Claim}
\newtheorem{proposition}[theorem]{Proposition}

\theoremstyle{definition}
\newtheorem{definition}[theorem]{Definition}

\newtheorem{remark}[theorem]{Remark}
\newtheorem{ex}[theorem]{Example}
\newtheorem{rem}[theorem]{Remark}

\numberwithin{equation}{section}

\graphicspath{{Images/}}

\newcommand{\EE}{\mathbb{E}}

\newcommand{\R}{\mathbb{R}}
\newcommand{\Z}{\mathbb{Z}}
\newcommand{\mcD}{\mathcal{D}}

\newcommand{\mcN}{\mathcal{N}}

\newcommand{\brr}{\overline}
\newcommand{\rb}{\mathbf{w}}

\renewcommand{\ll}{\llbracket}
\newcommand{\rr}{\rrbracket}
\newcommand{\mfh}{\mathcal{Q}_{2}^H}
\newcommand{\gj}{\mathcal{G}_H}
\newcommand{\gk}{\mathcal{G}_{r}}
\newcommand{\rev}{\textcolor{black}}
\newcommand{\bs}{\mathbf{s}}
\newcommand{\bw}{\mathbf{w}}
\newcommand{\s}{\mathbf{s}}

\newcommand{\mcL}{\mathcal{L}}
\newcommand{\vr}{\operatorname{Var}}
\newcommand{\Aut}{\operatorname{Aut}}
\newcommand{\Var}{\operatorname{Var}}
\newcommand{\Cov}{\operatorname{Cov}}
\newcommand{\Ex}{\mathbb{E}}
\newcommand{\ind}{\mathbf{1}}
\renewcommand{\P}{\mathbb{P}}
\newcommand{\skr}{\Upsilon_H}
\newcommand{\mth}{\mathcal{Q}_4^H}

\usepackage[normalem]{ulem}

\makeatletter
\providecommand*{\Dashv}{
  \mathrel{
    \mathpalette\@Dashv\vDash
  }
}
\makeatletter
\newcommand\thankssymb[1]{\textsuperscript{\@fnsymbol{#1}}}
\makeatother

\makeatother
\title[Asymptotic normality for monochromatic subgraphs]{A fourth moment phenomenon for asymptotic normality of monochromatic subgraphs
}
\author[S.\ Das]{Sayan Das\thankssymb{1}}
\address{S.\ Das, 
	Department of Mathematics, Columbia University
	}
\thanks{\thankssymb{1} Department of Mathematics, Columbia University, 2990 Broadway, New York, NY 10027, USA}
\email{sayan.das@columbia.edu}
\author[Z.\ Himwich]{Zoe Himwich\thankssymb{1}}
\address{Z.\ Himwich, 
	Department of Mathematics, Columbia University
	}
\email{himwich@math.columbia.edu}
\author[N.\ Mani]{Nitya Mani\thankssymb{2}}
\address{N.\ Mani, 
	Department of Mathematics, Massachusetts Institute of Technology
}
\thanks{\thankssymb{2} Department of Mathematics, Massachusetts Institute of Technology, 182 Memorial Drive, Cambridge, MA 02139, USA}
\email{nmani@mit.edu}

\begin{document}
\begin{abstract}
Given a graph sequence $\{G_n\}_{n\ge1}$ and a simple connected subgraph $H$, we denote by $T(H,G_n)$ the number of monochromatic copies of $H$ in
a uniformly random vertex coloring of $G_n$ with $c \ge 2$ colors. We prove a central limit theorem for $T(H,G_n)$ (we denote the appropriately centered and rescaled statistic as $Z(H,G_{n})$) with explicit error rates. The error rates arise from graph counts of collections formed by joining copies of $H$ which we call \textit{good joins}. Good joins are closely related to the fourth moment of $Z(H,G_{n})$, which allows us to show a \textit{fourth moment phenomenon} for the central limit theorem. For $c\ge 30$, we show that $Z(H,G_n)$ converges in distribution to $\mathcal{N}(0,1)$ whenever its fourth moment converges to $3$. We show the convergence of the fourth moment is necessary to obtain a normal limit when $c\ge2$.
\end{abstract}
\maketitle

\section{Introduction}

 Consider a deterministic sequence of simple undirected graphs $G_n= (V(G_n), E(G_n))$, where graph $G_n$ has vertex set $V(G_n) = \{1, 2, \ldots, v(G_n)\}$ and edge set $E(G_n)$. We denote the adjacency matrix of $G_{n}$ by $A(G_n) = (a_{ij}(G_n))_{i, j \in V(G_n)}$, where $a_{ij} = \ind\{(i, j) \in E(G_n)\}$.  We color the vertices of $G_{n}$ uniformly at random and independently with $c$ colors, so that
$$\P[X_v = a] = \frac{1}{c} \mbox{ for } a \in \{1,2,\ldots , c\},$$
where $X_v$ denotes the color of vertex $v \in V(G_n)$. Given a simple connected fixed subgraph $H$, we say that a copy of $H$ is \textit{monochromatic} if all of its vertices have the same color, i.e. $X_{u}=X_v$ for all $u,v\in V(H)$.  In this paper, we propose a general framework to study asymptotic normality for the statistic $T(H,G_n)$, the number of monochromatic copies of $H$ in $G_n$.

When $H = K_2$ is an edge, the statistic $T(K_2,G_n)$ counts the number of
monochromatic edges in $G_n$. 
This statistic arises in a variety of contexts. As one example, $T(K_2,G_n)$ counts the number of pairs of friends on friendship network $G_n$ that have the same birthday, thus generalizing the classical birthday paradox (which takes $G = K_n$). Consequently, the asymptotics of $T(K_2,G_{n})$ arise in various works studying generalizations of the birthday paradox~\cite{BHJ92,AGK16,CF06, D05,DM89} including the discrete logarithm problem~\cite{KMPT10,GH12} and problems related to testing discrete distributions~\cite{BFR13}. The asymptotics of $T(K_2, G_n)$ also naturally arise in a variety of applications, including understanding coincidences, cryptology, and occupancy urns~\cite{DM89,BHA16,NS07}.

The asymptotic distribution of $T(K_2,G_n)$ exhibits various universality phenomena, depending on how the number of colors $c$ scales with the number of edges $|E(G_n)|$. In particular, for fixed $c\ge 2$, using the method of moments, Bhattacharya, Diaconis, and Mukherjee~\cite{BHA17} concluded that $T(K_2,G_n)$ is asymptotically normal whenever $|E(G_n)|\to \infty$ and the fourth moment of a suitably normalized version of $T(K_2,G_n)$ converges to $3$ (the fourth moment of the normal distribution). The principle that a central limit theorem for nonlinear functions of random fields can arise as a consequence of convergence of the corresponding sequence of fourth moments arises in a variety of problems. This phenomenon was first observed in \cite{NP05,NP09}; we refer to Remark~\ref{remark:4moments} for a detailed discussion. Later, Bhattacharya, Fang, and Yan~\cite{BHA22} gave a limit theorem for the count of monochromatic triangles in family of graphs $G_n$, described by the statistic $T(K_3,G_n)$. By writing this statistic as martingale difference sequence, they were able to give a central limit theorem for $T(K_3,G_n)$, with an associated error bound for any fixed number of colors $c \ge 2$. When $c \ge 5$, these error bounds exhibit a similar fourth moment phenomenon to $T(K_2, G_n)$, but the authors observed that such a phenomenon does not necessarily arise for $2 \le c \le 4$. For more general $H$, Bhattacharya and Mukherjee~\cite{BHA19} obtained asymptotic results for $T(H,G_n)$, but only under the additional assumption that $G_n$ is dense.

In the Poisson regime, when $\Ex[T(H,G_n)]$ is bounded, the limiting distribution $T(H,G_n)$ is far better understood. Asymptotics have been shown for counts of edges~\cite{BHJ92,BHA17}, stars~\cite{star}, and for more general subgraphs~\cite{bmm20}. 
The corresponding questions in the Gaussian setting are less studied. Apart from the special cases of $H = K_2, K_3$, there has been little prior work giving precise conditions describing when $T(H,G_n)$ has a Gaussian limit for a general sequence of (not necessarily dense) graphs $G_n$.

In this paper, we initiate a systematic study of $T(H,G_n)$ for connected subgraphs $H$. The principal contributions of this article are:

\begin{enumerate}
    \item We prove a central limit theorem with an explicit error bound for the normalized version of $T(H,G_n)$, which we call $Z(H,G_n)$ (Theorem \ref{t:clt}). The error bounds hold for every $c\ge 2$, every sequence of deterministic graphs $G_n$ and for every connected subgraph $H$. 
    \item We provide a clear description of the graph counts which contribute to the error rates in this central limit theorem. We refer to these graph counts as \textit{good joins} and \textit{2-shared $2$-joins} of $H$. We introduce and discuss these definitions in Section~\ref{sec1.1}. 
    \item We show that the difference between the fourth moment of $Z(H,G_{n})$ and $3$ (the fourth moment of the normal distribution) is governed by counts of the same collection of subgraphs.
    This leads to a \textit{fourth moment phenomenon}. For $c\ge 30$, we show that $Z(H,G_n)$ is asymptotically normal whenever $\Ex[Z(H,G_n)^4]\to 3$ (Theorem~\ref{t:fourthmoment}). This confirms the prediction made by~\cite{BHA22} about fourth moment phenomena for general subgraphs. Interestingly, the lower bound `$30$' does not depend on $H$.
    \item We show that the fourth-moment condition is necessary to have a normal limit for any $c\ge 2$ fixed (Theorem~\ref{t:converse}). This yields asymptotic normality if and only if the fourth moment converges to $3$ whenever $c\ge 30$.
\end{enumerate}
  These results provide a comprehensive characterization of the asymptotic normality of $T(H,G_n)$ for any connected $H$.

\subsection{The model and the main results}\label{sec1.1}

Given a subgraph $H$, recall that $T(H,G_n)$ denotes the number of monochromatic copies of $H$ in $G_n$, 
\begin{align}\label{d:thgn}
T(H, G_n) := \frac{1}{|\Aut(H)|} \sum_{\bs \in V(G_n)_{|V(H)|}} \ind\{X_{=\bs}\} \prod_{(i, j) \in E(H)} a_{s_i s_j}(G_n),
\end{align}
where
\begin{itemize}
\setlength\itemsep{0.3em}
    \item[--] $V(G_n)_k= \{ \bs = (s_1, \ldots, s_{k}) \in V(G_n)^k : s_i \text{ distinct} \}$.
    \item[--] $\Aut(H)$ is the automorphism group of $H$, or the vertex permutations $\sigma$ such that $(i, j) \in E(H)$ if and only if $(\sigma(i), \sigma(j)) \in E(H)$.
    \item[--] $\ind\{X_{=\bs}\}:=\ind\{X_{s_1} = \cdots = X_{s_{k}}\}$ is the indicator variable that all of the vertices of $\bs\in V(G_n)_k$ have the same color.
\end{itemize}
For the rest of the paper, we assume $H$ is simple and connected and that $|V(H)|=r\ge 3$. Our main result provides central limit error rates in terms of \textit{joins} of graph counts.
\begin{definition}[$k$-join]
 A graph $F=(V(F),E(F))$ is a $k$-join of $H$ if there exist subgraphs $H_1,H_2,\ldots,H_k$ of $F$ such that each $H_i$ is isomorphic to $H$, $V(F)=V(\bigcup_{i=1}^k H_i)$, and $E(F)=E(\bigcup_{i=1}^k H_i)$.
 \end{definition}
Here we consider two special classes of joins of a given graph $H$:  \textit{good joins} and \textit{$2$-shared $2$-joins}.
\begin{definition}[Good join] \label{def:good}  Given a simple connected graph $H$, a \emph{good join} of $H$ is a graph formed by joining $4$ copies of $H$, say $H_1,H_2,H_3,H_4$, such that for all  $i=1,2,3,4$
		\begin{align} \label{good}
			\bigg| V(H_1\cup H_2\cup H_3 \cup H_4)\bigg|-\bigg| V\big(\bigcup_{j \neq i} H_j\big)\bigg| \le |V(H)|-2.
		\end{align}
		and
	\begin{align}\label{vgood}
		|V(H_1\cup H_2\cup H_3 \cup H_4)| \le \min_{\pi\in \mathbb{S}_4} \bigg\{ |V(H_{\pi_1}\cup H_{\pi_2})|+|V(H_{\pi_3}\cup H_{\pi_4})|-2 \bigg\}.
	\end{align}
	In the expression above, we take the minimum over all $\pi \in \mathbb{S}_4$, the set of permutations of $(1,2,3,4)$. We denote the collection of all good joins of $H$ by $\gj$.
\end{definition}
\begin{definition}[$2$-shared $k$-join] \label{2shared}   Given a simple, connected graph $H$, a \emph{$2$-shared $k$-join} of $H$ is a graph formed by joining $k$ copies of $H$, say $H_1,H_2,\ldots,H_k$,  such that $|V\left(\bigcap_{i=1}^k H_i\right)| \ge 2$. There are at least two distinct vertices which are present in all copies of $H_i$. We denote the collection of all $2$-shared $k$-joins of $H$ as $\mathcal{Q}_k^H$.
\end{definition}

The collections $\mathcal{Q}_k^H$ and $\gj$ both contain $H$. In either case, all the copies of $H$ in the join can completely overlap. Conditions~\eqref{good} and~\eqref{vgood} force a good join to be connected.  However as illustrated in Figures~\ref{fig:twostar}(A) and~\ref{fig:triangle}(A), a connected $4$-join is not necessarily a good join.
The $\mfh$ collection is a subcollection of $\gj$ (we can take $H_1=H_2$ and $H_3=H_4$ and arrange $H_1$ and $H_3$ with 2 common vertices). See Figures~\ref{fig:twostar} and~\ref{fig:triangle} for several examples of joins of two-stars $K_{1,2}$ and triangles $K_3$ (respectively).

Before we can state our results, we must introduce a couple more pieces of notation. We denote the centered and rescaled version of $T(H,G_n)$ by
\begin{align}\label{d:zhgn}
    Z(H, G_n) := \frac{T(H, G_n) - \EE[T(H, G_n)]}{\sigma_H(G_n)},\quad \sigma_H(G_n) := \sqrt{\Var[T(H, G_n)]}.
\end{align}
 Given a graph $F$, the number of copies of $F$ in $G_n$ is
\begin{align*}
    N(F,G_n):=\frac1{|\Aut(F)|}\sum_{s\in V(G_n)|_{|V(F)|}}\prod_{(i,j)\in E(F)} a_{s_is_j}(G_n).
\end{align*}
For a collection of graphs $\mathcal{H}$, $N(\mathcal{H},G_n):=\sum_{F\in \mathcal{H}} N(F,G_n)$ denotes the total number of copies of all graphs in $\mathcal{H}$ that can be found in $G_n$. We assume that $N(H,G_n) \ge 1$ for all $n\ge 1$. We are now ready to state our main result.
\begin{theorem}\label{t:clt} 
Fix a connected graph $H$ with $|V(H)|=r\ge 3$ and consider $Z(H,G_n)$ defined in \eqref{d:zhgn}. For all $c\ge 2$, there exist a constant $K= K(r,c)>0$ such that 
\begin{align}\label{e:clt}
\sup_{x\in \R} |\mathbb{P}[Z(H,G_n) \le x]-\Phi(x)| \le K\left(\frac{N(\gj,G_n)}{N(\mfh,G_n)^2}\right)^{1/20} 
\end{align}
where $\gj$ and $\mfh$ denote the collections of good joins and $2$-shared $2$-joins of $H$ defined in Definition~\ref{def:good} and~\ref{2shared} respectively. Here $\Phi$ denotes the cumulative distribution function of $\mcN(0,1)$, the standard normal distribution.
\end{theorem}

\begin{figure}[t]
    \centering
   \begin{subfigure}[b]{0.45\textwidth}
    \centering
    \begin{tikzpicture}[
    scale=0.75,
    font=\footnotesize,
    v/.style = {circle, fill, inner sep = 0pt, minimum size = 5pt},
    bv/.style = {circle, fill,blue, inner sep = 0pt, minimum size = 5pt},
    e/.style = {thick, black},
    be/.style = {thick, blue}]
    \node[label=below:$H_1$] (a) at (1,0) {};
    \node[label=below:$H_2$] (b) at (3,0) {};
    \node[label=below:$H_3$] (c) at (5,0) {};
    \node[label=below:$H_4$] (d) at (7,0) {};
    \node[v,label=below:1] (1) at (0,0) {};
    \node[v,label=right:2] (2) at (1,1.73) {};
    \node[v,label=below:3] (3) at (2,0) {};
    \node[v,label=right:4] (4) at (3,1.73) {};
    \node[v,label=below:5] (5) at (4,0) {};
    \node[v,label=right:6] (6) at (5,1.73) {};
    \node[v,label=below:7] (7) at (6,0) {};
    \node[bv,label=right:8] (8) at (7,1.73) {};
    \node[bv,label=below:9] (9) at (8,0) {};
    \draw[e] (1)--(2)--(3)--(4)--(5)--(6)--(7)--(8);
    \draw[be] (8)--(9);
    \end{tikzpicture}
    \caption{}
   \end{subfigure}
   \begin{subfigure}[b]{0.45\textwidth}
   \centering
   \begin{tikzpicture}[
    scale=0.75,
    font=\footnotesize,
    v/.style = {circle, fill, inner sep = 0pt, minimum size = 5pt},
    bv/.style = {circle, fill,blue, inner sep = 0pt, minimum size = 5pt},
    e/.style = {thick, black},
    be/.style = {thick, blue}]
    \node[label=above left:$H_1$] (a) at (0.5,0.8) {};
    \node[label=below left:$H_2$] (b) at (0.5,-0.8) {};
    \node[label=above right:$H_3$] (c) at (3.5,0.8) {};
    \node[label=below right:$H_4$] (d) at (3.5,-0.8) {};
    \node[v,label=below:1] (1) at (0,0) {};
    \node[v,label=right:2] (2) at (1,1.73) {};
    \node[bv,label=below:3] (3) at (2,0) {};
    \node[v,label=right:4] (4) at (3,1.73) {};
    \node[v,label=below:5] (5) at (4,0) {};
    \node[v,label=right:6] (6) at (1,-1.73) {};
    \node[v,label=right:7] (7) at (3,-1.73) {};
    \draw[e] (6)--(1)--(2)--(3)--(4)--(5)--(7)--(3)--(6);
    \end{tikzpicture}
    \caption{}
   \end{subfigure}
      \begin{subfigure}[b]{0.3\textwidth}
    \centering
 \begin{tikzpicture}[
    scale=0.75,
    font=\footnotesize,
    v/.style = {circle, fill, inner sep = 0pt, minimum size = 5pt},
    bv/.style = {circle, fill,red, inner sep = 0pt, minimum size = 5pt},
    e/.style = {thick, black},
    be/.style = {thick, red}]
    \node[label=below:$H_1$] (a) at (1,0) {};
    \node[label=below:$H_2$] (b) at (3,0) {};
    \node[label=above:$H_3$] (c) at (2,1.73) {};
    \node[label=above:$H_4$] (d) at (4,1.73) {};
    \node[v,label=below:1] (1) at (0,0) {};
    \node[v,label=above:2] (2) at (1,1.73) {};
    \node[v,label=below:3] (3) at (2,0) {};
    \node[v,label=above:4] (4) at (3,1.73) {};
    \node[v,label=below:5] (5) at (4,0) {};
    \node[v,label=above:6] (6) at (5,1.73) {};
    \draw[e] (1)--(2)--(3)--(4)--(5)--(6);
    \end{tikzpicture}
    \caption{}
   \end{subfigure}
   \begin{subfigure}[b]{0.3\textwidth}
   \centering
   \begin{tikzpicture}[
    scale=0.75,
    font=\footnotesize,
    v/.style = {circle, fill, inner sep = 0pt, minimum size = 5pt},
    bv/.style = {circle, fill,red, inner sep = 0pt, minimum size = 5pt},
    e/.style = {thick, black},
    be/.style = {thick, red}]
    \node[label=above left:$H_1$] (a) at (0.4,0.7) {};
    \node[label=below left:$H_2$] (b) at (3,-0.7) {};
    \node[label=above:$H_3$] (c) at (1,0) {};
    \node[label=below:$H_4$] (d) at (1,-0.) {};
    \node[v,label=left:1] (1) at (0,0) {};
    \node[v,label=above:2] (2) at (1,1) {};
    \node[v,label=right:3] (3) at (2,0) {};
    \node[v,label=below:4] (4) at (-1.5,0) {};
    \node[v,label=below:5] (5) at (3.5,0) {};
    \node[v,label=below:6] (6) at (1,-1) {};
   \draw[e] (6)--(1)--(2)--(4)--(6)--(5)--(2)--(3)--(6);
    \end{tikzpicture}
    \caption{}
   \end{subfigure}
    \begin{subfigure}[b]{0.3\textwidth}
   \centering
   \begin{tikzpicture}[
    scale=0.75,
    font=\footnotesize,
    v/.style = {circle, fill, inner sep = 0pt, minimum size = 5pt},
    bv/.style = {circle, fill,red, inner sep = 0pt, minimum size = 5pt},
    e/.style = {thick, black},
    be/.style = {thick, red}]
    \node[label=above:$H_1$] (c) at (1,0) {};
    \node[label=below:$H_2$] (d) at (1,-0.) {};
    \node[v,label=left:1] (1) at (0,0) {};
    \node[v,label=above:2] (2) at (1,1) {};
    \node[v,label=right:3] (3) at (2,0) {};
    \node[v,label=below:4] (4) at (1,-1) {};
   \draw[e] (1)--(2)--(3)--(4)--(1);
    \end{tikzpicture}
    \caption{}
   \end{subfigure}
    \caption{Examples of joins for $2$-star $K_{1,2}$: (A) is a connected $4$-join of $K_{1,2}$ which does not satisfy \eqref{good}. (B) is a connected $4$-join of $K_{1,2}$ that satisfies~\eqref{good} but not~\eqref{vgood}. (C) and (D) are examples of good joins of $K_{1,2}$. (E) is an example of $2$-shared $2$-join of $K_{1,2}$.}
    \label{fig:twostar}
\end{figure}

\begin{figure}[t]
    \centering
   \begin{subfigure}[b]{0.45\textwidth}
    \centering
    \begin{tikzpicture}[
    scale=0.75,
    font=\footnotesize,
    v/.style = {circle, fill, inner sep = 0pt, minimum size = 5pt},
    bv/.style = {circle, fill,blue, inner sep = 0pt, minimum size = 5pt},
    e/.style = {thick, black},
    be/.style = {thick, blue}]
    \node[label=below:$H_1$] (a) at (1,0) {};
    \node[label=below:$H_2$] (b) at (3,0) {};
    \node[label=below:$H_3$] (c) at (5,0) {};
    \node[label=below:$H_4$] (d) at (7,0) {};
    \node[v,label=below:1] (1) at (0,0) {};
    \node[v,label=right:2] (2) at (1,1.73) {};
    \node[v,label=below:3] (3) at (2,0) {};
    \node[v,label=right:4] (4) at (3,1.73) {};
    \node[v,label=below:5] (5) at (4,0) {};
    \node[v,label=right:6] (6) at (5,1.73) {};
    \node[v,label=below:7] (7) at (6,0) {};
    \node[bv,label=right:8] (8) at (7,1.73) {};
    \node[bv,label=below:9] (9) at (8,0) {};
    \draw[e] (1)--(2)--(3)--(4)--(5)--(6)--(7)--(8);
    \draw[be] (8)--(9);
    \draw[e] (1)--(3)--(5)--(7)--(9);
    \end{tikzpicture}
    \caption{}
   \end{subfigure}
   \begin{subfigure}[b]{0.45\textwidth}
   \centering
   \begin{tikzpicture}[
    scale=0.75,
    font=\footnotesize,
    v/.style = {circle, fill, inner sep = 0pt, minimum size = 5pt},
    bv/.style = {circle, fill,blue, inner sep = 0pt, minimum size = 5pt},
    e/.style = {thick, black},
    be/.style = {thick, blue}]
    \node[label=above left:$H_1$] (a) at (0.5,0.8) {};
    \node[label=below left:$H_2$] (b) at (0.5,-0.8) {};
    \node[label=above right:$H_3$] (c) at (3.5,0.8) {};
    \node[label=below right:$H_4$] (d) at (3.5,-0.8) {};
    \node[v,label=below:1] (1) at (0,0) {};
    \node[v,label=right:2] (2) at (1,1.73) {};
    \node[bv,label=below:3] (3) at (2,0) {};
    \node[v,label=right:4] (4) at (3,1.73) {};
    \node[v,label=below:5] (5) at (4,0) {};
    \node[v,label=right:6] (6) at (1,-1.73) {};
    \node[v,label=right:7] (7) at (3,-1.73) {};
    \draw[e] (6)--(1)--(2)--(3)--(4)--(5)--(7)--(3)--(6);
    \draw[e] (1)--(3)--(5);
    \end{tikzpicture}
    \caption{}
   \end{subfigure}
      \begin{subfigure}[b]{0.3\textwidth}
    \centering
 \begin{tikzpicture}[
    scale=0.75,
    font=\footnotesize,
    v/.style = {circle, fill, inner sep = 0pt, minimum size = 5pt},
    bv/.style = {circle, fill,red, inner sep = 0pt, minimum size = 5pt},
    e/.style = {thick, black},
    be/.style = {thick, red}]
    \node[label=below:$H_1$] (a) at (1,0) {};
    \node[label=below:$H_2$] (b) at (3,0) {};
    \node[label=above:$H_3$] (c) at (2,1.73) {};
    \node[label=above:$H_4$] (d) at (4,1.73) {};
    \node[v,label=below:1] (1) at (0,0) {};
    \node[v,label=above:2] (2) at (1,1.73) {};
    \node[v,label=below:3] (3) at (2,0) {};
    \node[v,label=above:4] (4) at (3,1.73) {};
    \node[v,label=below:5] (5) at (4,0) {};
    \node[v,label=above:6] (6) at (5,1.73) {};
    \draw[e] (1)--(2)--(3)--(4)--(5)--(6);
    \draw[e] (1)--(3)--(5);
    \draw[e] (2)--(4)--(6);
    \end{tikzpicture}
    \caption{}
   \end{subfigure}
   \begin{subfigure}[b]{0.3\textwidth}
   \centering
   \begin{tikzpicture}[
    scale=0.75,
    font=\footnotesize,
    v/.style = {circle, fill, inner sep = 0pt, minimum size = 5pt},
    bv/.style = {circle, fill,red, inner sep = 0pt, minimum size = 5pt},
    e/.style = {thick, black},
    be/.style = {thick, red}]
    \node[label=above left:$H_1$] (a) at (0.4,0.7) {};
    \node[label=below left:$H_2$] (b) at (3,-0.7) {};
    \node[label=above:$H_3$] (c) at (0.5,-0.5) {};
    \node[label=below:$H_4$] (d) at (1.5,0.5) {};
    \node[v,label=left:1] (1) at (-0.5,0) {};
    \node[v,label=above:2] (2) at (1,1) {};
    \node[v,label=right:3] (3) at (2.5,0) {};
    \node[v,label=below:4] (4) at (-1.5,0) {};
    \node[v,label=below:5] (5) at (3.5,0) {};
    \node[v,label=below:6] (6) at (1,-1) {};
   \draw[e] (6)--(1)--(2)--(4)--(6)--(5)--(2)--(3)--(6);
   \draw[e] (2)--(6);
    \end{tikzpicture}
    \caption{}
   \end{subfigure}
    \begin{subfigure}[b]{0.3\textwidth}
   \centering
   \begin{tikzpicture}[
    scale=0.75,
    font=\footnotesize,
    v/.style = {circle, fill, inner sep = 0pt, minimum size = 5pt},
    bv/.style = {circle, fill,red, inner sep = 0pt, minimum size = 5pt},
    e/.style = {thick, black},
    be/.style = {thick, red}]
    \node[label=above:$H_1$] (c) at (1,0) {};
    \node[label=below:$H_2$] (d) at (1,-0.) {};
    \node[v,label=left:1] (1) at (0,0) {};
    \node[v,label=above:2] (2) at (1,1.73) {};
    \node[v,label=right:3] (3) at (2,0) {};
    \node[v,label=below:4] (4) at (1,-1.73) {};
   \draw[e] (1)--(2)--(3)--(4)--(1)--(3);
    \end{tikzpicture}
    \caption{}
   \end{subfigure}
    \caption{Examples of joins for the triangle $K_{3}$: (A) is a connected $4$-join of $K_{3}$ which does not satisfy \eqref{good}. (B) is a connected $4$-join of $K_{3}$ that satisfies~\eqref{good} but not~\eqref{vgood}. (C) and (D) are examples of good joins of $K_{3}$. (E) is an example of $2$-shared $2$-join of $K_{3}$.}
    \label{fig:triangle}
\end{figure}

Theorem~\ref{t:clt} implies that $Z(H,G_n) \overset{\mathcal{D}}\implies \mcN(0,1)$ whenever $$\frac{N(\gj, G_n)}{N(\mfh,G_n)^2}\to 0.$$ The proof of Theorem \ref{t:clt} proceeds by writing $Z(H,G_n)$ as a martingale difference sequence. This decomposition result (Proposition~\ref{p:zrgn0}) is the starting point of our analysis. Proposition~\ref{p:zrgn0} can be derived from the Hoeffding decomposition~\cite{HOE48}. Instead, we directly prove Proposition~\ref{p:zrgn0} by verifying martingale properties of the decomposition in Section~\ref{sec2}. By applying the standard martingale CLT with error bounds, a result in~\cite{HEY70}, we obtain error rates in terms of moments of the martingale sequence. One of the main contributions of this paper is a delicate understanding of these error rates in terms of graph counts. We show that the nontrivial contributions from the errors precisely come from the good joins of $H$ introduced in Definition~\ref{def:good}.

Our error rates are also related to the fourth moment of $Z(H,G_n)$ which leads to the following result.

\begin{theorem}\label{t:fourthmoment}
Fix a connected graph $H$ with $|V(H)|=r\ge 3$ and consider $Z(H,G_n)$ defined in \eqref{d:zhgn}. For all $c\ge 30$, there exist a constant $K= K(r,c)>0$ such that  
\begin{align*}
\sup_{x\in \R} |\mathbb{P}[Z(H,G_n) \le x]-\Phi(x)| \le K\left(\Ex[Z(H,G_n)^4]-3\right)^{\frac1{20}},
\end{align*}
where $\Phi$ denotes the cumulative distribution function of $\mcN(0,1)$.
\end{theorem}

The above result shows that the asymptotic normality of normalized counts of monochromatic copies of $H$ are governed by a \textit{fourth moment phenomenon} (see Remark \ref{remark:4moments} for a discussion). Indeed, Theorem \ref{t:fourthmoment} claims that for $c\ge 30$, when $\Ex[Z(H,G_n)^4]$ converges to $3$ (the fourth moment of normal distribution), $Z(H,G_n)$ converges weakly to the standard normal distribution. One caveat of the above theorem is that it requires $c\ge 30$. A nontrivial lower bound on the number of colors is necessary. Indeed, it was demonstrated in~\cite[Theorem 2~(2)]{BHA22} that in the case of triangles, for $2\le c\le 4$,  fourth moment convergence to $3$ may not be sufficient for asymptotic normality. 

\begin{rem}
The lower bound `$30$' for $c$ given in Theorem~\ref{t:fourthmoment} is a uniform lower bound for any connected graph $H$. One interesting direction for future work could be determining the optimal lower bound for a specific $H$. Further, we study $T(H, G_n)$ only in the regime where $c$ is fixed. In the setting where $c = c_n\to \infty$, \cite{BHA17,xiao} gave asymptotic normality results for $H=K_2$. Another interesting direction for future work could be to extend the results in this paper to the asymptotic case where $c_{n}\to\infty$. 
\end{rem}

The proof of Theorem~\ref{t:fourthmoment} shows that the good join counts $N(\gj,G_n)$ are exactly the terms which contribute to $\Ex[Z(H,G_n)^4]-3$. For triangles, this result appears as~\cite[Theorem 2~(1)]{BHA22}, with $c\ge 5$. In~\cite{BHA22}, the authors provide exact computations of $\Ex[Z({K}_3,G_n)^4]-3$ and observe $32$ different $4$-joins of ${K}_3$ which contribute. It is not hard to check that their collection of $4$-joins is contained in $\mathcal{G}_{{K}_3}$. In the case of a general $H$, we provide a careful understanding of mixed central moments of the color indicators $\ind\{X_{=\bs}\}$ (see Lemma~\ref{l:zexpec} and Lemma~\ref{l:not-very-good-factors}). At the expense of a larger $c$, these lead to a rather elegant proof of Theorem~\ref{t:fourthmoment} that works for all $H$.

A natural follow up question is whether the converse of Theorem \ref{t:fourthmoment} holds. Our next result answers this question affirmatively.

\begin{theorem}\label{t:converse}
For $c\ge 2$ fixed, we have that $Z(H,G_n) \overset{\mcD}\implies \mcN(0, 1)$ implies $\EE[Z(H,G_n)^4] \rightarrow 3$.
\end{theorem}

This result, along with Theorem \ref{t:fourthmoment}, shows that for all $c\ge 30$, asymptotic normality for $Z(H,G_n)$ occurs if and only if the fourth moment of $Z(H,G_n)$ converges to $3$. To prove Theorem~\ref{t:converse}, we find that all higher moments of $Z(H,G_n)$ are uniformly bounded in $n$ (Theorem \ref{thm:mom-bd}). Theorem~\ref{t:converse} then follows by an application of uniform integrability. To prove Theorem~\ref{thm:mom-bd}, we adopt the approach of~\cite{BHA22}, which relies on estimates from extremal combinatorics.

\begin{remark}\label{remark:4moments}
The fourth moment phenomenon was first observed for Wiener-It\^o stochastic integrals in~\cite{NP05}, and~\cite{NP09} later provided error bounds for the fourth moment result in~\cite{NP05}. 
In subsequent years, the fourth moment approach has appeared as a governing principle for many central limit results for non-linear functionals of random fields.  The book~\cite{NP12} provides a survey on the topic and the website \url{https://sites.google.com/site/malliavinstein/home} holds an up-to-date list of related results. 
In the context of coloring problems, fourth moment phenomena have been observed for monochromatic edges~\cite{BHA17} and for triangles~\cite{BHA22}. Instead of considering monochromatic copies of $H$, one can also consider monochromatic copies of $H$ of a fixed color. In this framework,~\cite{bdm22} observed that a certain fourth moment phenomenon holds for counts of monochromatic copies of $H$ of a fixed color. 
 \end{remark}

\subsection{Assumptions and Notation} 
We write $v(F)$ to denote the number of vertices of a graph $F$. For the remainder of the paper we work with a deterministic sequence of graphs $G_n$, where $G_n=(V(G_n),E(G_n))$ with adjacency matrix $A(G_n) = (a_{ij}) \in \{0, 1\}^{v(G_n) \times v(G_n)}$. We work with a fixed connected subgraph $H$ with $v(H)=r\ge 3$ and a fixed number of colors $c\ge 2$. We assume $N(H,G_n)\ge 1$ for all $n\ge 1$. For $\bs \in V(G_n)_{k}$, we set $\brr\s=\{s_1,s_2,\ldots,s_k\}$ and define $$a_{H, \bs} := \prod_{(i, j) \in E(H)} a_{s_is_j}.$$
Note that $a_{H,\bs}$ depends on $G_n$. We will often hide this dependency. 

\smallskip

Throughout the paper, we use standard asymptotic notation. In particular, $f(n) \lesssim_{\square} g(n)$ and $f(n) \gtrsim_{\square} g(n)$ denote that $f(n) \le C_1 \cdot g(n)$  and $f(n)\ge C_2\cdot g(n)$ for $C_1,C_2 > 0$ which depend on the subscript parameters. We write $f(n)\asymp_{\square} g(n)$ if $f(n)\lesssim_{\square} g(n)$ and $f(n) \gtrsim_{\square} g(n)$ holds simultaneously. It is important to note that the $\lesssim$ and $\gtrsim$ appearing in later sections will always depend on $r$ and $c$. We won't mention this further. We will only highlight the dependencies on other parameters. Finally, for $a,b\in \R$ we use $\ll a,b\rr:=[a,b]\cap \Z$, $a \land b = \min(a, b)$, and $[a] = \ll 1,a \rr$. 

\subsection*{Organization} In Section~\ref{sec2} we provide a decomposition of $Z(H,G_n)$ as a martingale difference sequence. In Section \ref{sec3} we prove Theorem~\ref{t:clt}. Proofs of Theorem~\ref{t:fourthmoment} and Theorem~\ref{t:converse} are given in Section~\ref{sec4} and Section~\ref{sec5} respectively.

\subsection*{Acknowledgements} We thank Sumit Mukherjee for suggesting this problem. We thank Bhaswar Bhattacharya and Sumit Mukherjee for useful discussions. \rev{We thank the anonymous referees for their careful reading and useful comments on improving our
manuscript.} SD acknowledges support from NSF DMS-1928930 during his participation in
the program ``Universality and Integrability in Random Matrix Theory and Interacting Particle Systems'' hosted
by the Mathematical Sciences Research Institute in Berkeley, California in fall 2021.  SD’s research was also partially supported by the Fernholz Foundation’s ``Summer Minerva Fellows'' program. NM was supported by a Hertz Graduate Fellowship and the NSF Graduate Research Fellowship Program. 

\section{A decomposition for $Z(H,G_n)$} \label{sec2}

 The goal of this section  is to prove Proposition~\ref{p:zrgn0}, which gives a decomposition of the normalized subgraph counts $Z(H,G_n)$ (defined in~\eqref{d:zhgn}). To state the decomposition precisely, we need a few definitions. First, we introduce the notion of \textit{strictly increasing tuples} of vertices (as opposed to the distinct tuples $V(G_n)_k$ considered in~\eqref{d:thgn}) along with some other helpful notation.

\begin{definition}[Strictly increasing tuples] \label{d:skt} For $k\in \ll 2,r\rr$,
$$\Lambda_{k,*}:=\{\rb:=(w_1,w_2,\ldots,w_k)\in V(G_n)_k \mid 1\le w_1 < w_2 < \cdots <w_k\le v(G_n)\},$$ 
denotes collection of all strictly increasing $k$-tuples of vertices from $V(G_n)$. Since elements of $\rb \in \Lambda_{k,*}$ are always arranged in strictly increasing order, we can also view $\rb$ as a set.

For $t\in V(G_n)$ and $k\in \ll 2,r\rr$,
\begin{align}\label{e:lkt}
    \Lambda_{k,t}:=\{\rb:=(w_1,w_2,\ldots,w_k)\in \Lambda_{k,*} \mid w_k=t\}
\end{align} denotes the collection of strictly increasing $k$-tuples with last vertex fixed to be $t\in V(G_n)$. For $\rb_1,\rb_2\in \Lambda_{k,t}$ we write $\rb_1\neq \rb_2$ if they do not define the same set.

For any $k_1,k_2\in \ll 2,r \rr$, $t\in V(G_n)$ and $u \in \ll 1, k_1\wedge k_2\rr$, we set
\begin{align}\label{e:ovlap}
\Lambda_{k_1,k_2,t}^{(u)}:=\left\{ (\rb_1,\rb_2) \mid \rb_i \in \Lambda_{k_i,t}, |\rb_1 \cap \rb_2|=u\right\}.
\end{align}
This means that $\Lambda_{k_1,k_2,t}^{(u)}$ contains pairs of strictly increasing tuples of size $k_1$ and $k_2$ which have the same last vertex, $t$, and which share exactly $u$ elements. Note that $u\ge 1$ as $t$ is a common element.
\end{definition}
\begin{rem}
We remark that $\Lambda_{k,*}$ and $\Lambda_{k,t}$ both depend on $n$ via $V(G_n)$, but we have suppressed this dependence. We will typically reserve the letter $\rb$ to denote $k$-tuples from $\Lambda_{k,*}$ or $\Lambda_{k,t}$ (as opposed to $\bs$ which will be used for elements from $V(G_n)_k$). The value of $k$ or $t$ will always be clear from the context and therefore we suppress it from the notation $\rb$.
\end{rem}

Given integer $k \in [r]$ and $\rb=(w_1,\ldots,w_k) \in \Lambda_{k,*}$, let 
$D_{\rb}(G_n)$ be the number of copies of $H$ in $G_n$ that include vertices $w_1,\ldots,w_k$.
More precisely,
\begin{align}\label{d:dw}
D_{\rb}:=\frac1{\Aut(H)}\sum_{\s\in V(G_n)_r  :  \brr{\bs} \supset \rb}  a_{H,\s}.
\end{align}
Next, we define the building blocks of our decomposition of $Z(H,G_n):$ the random variables $Y_{\rb}$. For $\rb\in \Lambda_{k,*}$ we let 
\begin{align}\label{deftY}
    \widetilde{Y}_{\rb}:=D_{\rb} \cdot \sum_{p=2}^{k}(-1)^{k-p}\sum_{1\le i_1 < \cdots <i_p\le k}\frac1{c^{r-p}}\ind\{X_{w_{i_1}}=\cdots=X_{w_{i_p}}\}.
\end{align}
We also define a centered analogue:
\begin{align}\label{e:defY}
    {Y}_{\rb}:=D_{\rb} \cdot \sum_{p=2}^{k}(-1)^{k-p}\sum_{1\le i_1 < \cdots <i_p\le k}\left(\frac1{c^{r-p}}\ind\{X_{w_{i_1}}=\cdots=X_{w_{i_p}}\}-\frac{1}{c^{r-1}}\right).
\end{align}
Finally, for any $t \in V(G_n)$, we set 
\begin{align}\label{d:ut}
    U_{t}:=\frac{1}{\sigma_H}\sum_{k=2}^r \sum_{\rb \in \Lambda_{k,t}}Y_{\rb},
\end{align}
where $\Lambda_{k,t}$ is defined in \eqref{e:lkt} and $\sigma_H$ is defined in \eqref{d:zhgn}. 

\medskip

With the above notation in place, we state a decomposition result for $Z(H,G_n)$ (defined in \eqref{d:zhgn}).

\begin{proposition}[$Z(H,G_n)$ decomposition]\label{p:zrgn0} The random variables $U_t$ form a martingale difference sequence with respect to the filtration $\mathcal{F}_t$ generated by $X_1,X_2,\ldots,X_t$. In other words, $U_t$ is measurable with respect to $X_1,X_2,\ldots,X_t$ and $\Ex[U_{t}\mid X_1,X_2,\ldots,X_{t-1}]=0$. Furthermore, 
\begin{align}\label{ob0}
    Z(H,G_n)=\sum_{t=1}^{v(G_n)} U_t.
\end{align}
\end{proposition}

\begin{proof}[Proof of Proposition~\ref{p:zrgn0}]  Recall $U_t$ from \eqref{d:ut}. From its definition, it is clear that $U_t$ is a measurable function of $X_{1},\ldots,X_{t}$. To show that $U_t$ is a martingale difference sequence we rely on the following lemma.

\begin{lemma}\label{l:ycond} Fix any $k\in \ll 2,r\rr$. Recall $\Lambda_{k,*}$ from Definition~\ref{d:skt}. \rev{Fix any $\rb=(w_1,w_2,\ldots,w_k)\in \Lambda_{k,*}$. Let $\Gamma \subset V(G_n)$ such that ${\rb} \not\subset \Gamma$. We have
\begin{align*}
    \Ex[Y_{\rb}\mid X_v, v\in \Gamma]=0.
\end{align*}}
\end{lemma}
Note that $U_{t}$ is a linear combination of terms of the form $Y_{w_1 \ldots w_{k-1}t}$ for $1 \le w_1 < \cdots < w_{k-1} \le t-1$. Thus, the fact that $U_t$ is a martingale with zero expectation is immediate from Lemma~\ref{l:ycond}. To check the last part of Proposition~\ref{p:zrgn0}, that the $U_{t}$ form a decomposition of $Z(H,G_{n})$ via~\eqref{ob0}, we use the following result.

\begin{lemma}\label{l:zrgn} Recall $\Lambda_{k,*}$ from Definition~\ref{d:skt}. Then,
\begin{align}\label{ob1}
    T(H,G_n)=\sum_{k=2}^r\sum_{\rb\in \Lambda_{k,*}} \widetilde{Y}_{\rb},
\end{align}
where $\widetilde{Y}_{\rb}$ is defined in \eqref{deftY}. 
\end{lemma}
 Since $\EE[\ind\{X_{w_{i_1}}=\cdots=X_{w_{i_p}}\}] = \frac{1}{c^{p-1}}$ for $1 \le i_1 < \cdots < i_p \le k$, we subtract this expectation from both sides of \eqref{ob1} and divide by $\sigma_H$ to get 
$$Z(H,G_n)=\frac1{\sigma_H}\sum_{k=2}^r\sum_{\rb\in \Lambda_{k,*}} {Y}_{\rb}.$$
Fixing the largest vertex $w_k$ to be $t$ for each $k$, and then summing over other vertices first, we obtain the alternative representation of $Z(H,G_n)$ in \eqref{ob0}. This completes the proof of Proposition~\ref{p:zrgn0}, up to proving Lemmas~\ref{l:ycond} and~\ref{l:zrgn}.
\end{proof}
\begin{proof}[Proof of Lemma~\ref{l:ycond}]

Fix any $\rb=(w_1,w_2,\ldots,w_k)\in \Lambda_{k,*}$. \rev{Note that $Y_{\rb}$ defined in~\eqref{e:defY} is a symmetric function of $w_1,\ldots,w_k$ only. Hence to prove the lemma it suffices to show
\begin{align*}
    \Ex[Y_{\rb}\mid X_{w_1},\ldots,X_{w_{k-1}}]=0.
\end{align*}}
We first break up the terms in the definition of $Y_{\rb}$ based on whether or not they contain $w_k$.
Towards this end, we define $M_1$ to consist of the terms in the expansion of  $Y_{\rb}$ that do not contain $w_k$ (so we require $i_p < k$, $p < k$),
\begin{align*}
M_1 & := \sum_{p = 2}^{k-1} (-1)^{k-p} \sum_{1 \le i_1 < \cdots < i_p < k} \left(\frac{1}{c^{r-p}} \ind\{X_{w_{i_1}} = \cdots = X_{w_{i_p}}\} - \frac{1}{c^{r-1}}\right)
\end{align*}
$M_2$ comprises the terms in the expansion in \eqref{e:defY} that contain $w_k$, but with $p > 2$: 
\begin{align*}
M_2 :&= \sum_{p = 3}^{k} (-1)^{k-p} \sum_{1 \le i_1 < \cdots < i_{p-1} <i_p=k} \left( \frac{1}{c^{r-p}} \ind\{X_{w_{i_1}} = \cdots = X_{w_{i_{p-1}}}=X_{w_k}\} - \frac{1}{c^{r-1}} \right) 
\end{align*}
Lastly, the collection of terms which remain form $M_3$:
\begin{align*}
M_3 :&= (-1)^{k-2} \sum_{i = 1}^{k-1} \left(\frac{1}{c^{r-2}} \ind \{X_{w_i} = X_{w_k} \} - \frac{1}{c^{r-1}} \right). 
\end{align*}
Therefore,
$$Y_{\rb}=D_{\rb} \cdot (M_1+M_2+M_3).$$

We take the conditional expectation of the $M_i$, conditioning on $X_{w_1}, \ldots, X_{w_{k-1}}$. Observe that $M_1$ is measurable with respect to~ $X_{w_1}, \ldots, X_{w_{k-1}}$. For $M_2$ note that for any $1\le i_1\le \cdots \le i_{p-1}<k$
\begin{align*}
    \Ex[ \ind\{X_{w_{i_1}}=\cdots=X_{w_{i_{p-1}}}=X_{w_{k}}\} ] = \frac1c\ind\{X_{w_{i_1}}=\cdots=X_{w_{i_{p-1}}}\}.
\end{align*}
This implies
\begin{align*}
\Ex[M_2\mid X_{w_1},X_{w_2},\ldots,X_{w_{k-1}}]=\sum_{p = 3}^{k} (-1)^{k-p} \sum_{1 \le i_1 < \cdots < i_{p-1} <k} \left[ \frac{1}{c^{r-p+1}} \ind\{X_{w_{i_1}} = \cdots = X_{w_{i_{p-1}}}\} - \frac{1}{c^{r-1}} \right].
\end{align*}
Making the change of variable $q := p - 1$ we see that the above sum is the same as $-M_1$. Finally for $M_3$, note that $\P(X_{w_i}=X_{w_k}\mid X_{w_i})=\frac1c$ for all $i \in [k-1]$, forcing $\Ex[M_3\mid X_{w_1},\ldots,X_{w_{k-1}}]=0.$ Thus, as desired, $$\Ex[Y_{\rb}\mid X_{w_1},\ldots,X_{w_{k-1}}]=D_{\rb} \cdot \Ex[M_1+M_2+M_3\mid X_{w_1},\ldots,X_{w_{k-1}}]=0.$$
\end{proof}

\begin{proof}[Proof of Lemma~\ref{l:zrgn}] Recall $\widetilde{Y}_\rb$ from \eqref{deftY}. Using the definition of $D_{\rb}$, we have
$$\widetilde{Y}_{\rb} = \frac1{|\Aut(H)|}\sum_{\s\in V(G_n)_r} a_{H,\s} \ind\{\rb\subset \brr\s\} \cdot \widetilde{M}_{\rb},$$ where 
$$\widetilde{M}_{\rb}:=\sum_{p=2}^{k}(-1)^{k-p}\sum_{1\le i_1 < \cdots <i_p\le k}\frac1{c^{r-p}}\ind\{X_{w_{i_1}}=\cdots=X_{w_{i_p}}\}.$$
Interchanging the ordering of the sum on the right hand side of~\eqref{ob1}, we have
\begin{align*}
    \sum_{k=2}^r\sum_{\rb\in \Lambda_{k,*}} \widetilde{Y}_{\rb}= \frac1{\Aut(H)}\sum_{\s\in V(G_n)_r} a_{H,\s} \sum_{k=2}^r\sum_{\rb\in \Lambda_{k,*}} \ind\{\rb\subset \brr\s\} \cdot \widetilde{M}_{\rb}.
\end{align*}
Fix any $\s\in V(G_n)_r$. Choose any $p\in \ll 2,r-1\rr$ (note that $p$ is strictly less than $r$). Fix a set of $p$ distinct vertices $\{v_1, \ldots, v_p\}$. Let $C_{v_1, \ldots, v_p}(\bs)$ be the coefficient of $\ind\{X_{v_1}=X_{v_2}=\ldots =X_{v_p}\}$ appearing in the double sum 
\begin{align*}
    \sum_{k=2}^r\sum_{\rb\in \Lambda_{k,*}} \ind\{\rb\subset \brr\s\} \cdot \widetilde{M}_{\rb}.
\end{align*}
We claim that $C_{v_1, \ldots, v_p}(\bs)= 0$. Clearly given any $\rb \in \Lambda_{k,*}$, the indicator appears in $\widetilde{M}_{\rb}$ if and only if $\{v_1,v_2,\ldots,v_p\} \subset \rb$, doing so with coefficient $(-1)^{k-p}c^{p-r}$.  Consequently, we observe that
\begin{align*}
C_{v_1, \ldots, v_p}(\bs) &= \sum_{k = 2}^r \sum_{\rb\in \Lambda_{k,*}} \ind\{\{v_1,v_2,\ldots,v_p\} \subset \rb \subset \brr\s \} (-1)^{k-p} c^{p-r}. 
\end{align*}
Clearly $\ind\{\{v_1,v_2,\ldots,v_p\} \subset \rb \subset \brr\s \}  = 0$ for $\rb\in \Lambda_{k,*}$ with $k<p$. For $k\ge p$, given any $\bs\in V(G_n)_r$ and $\{v_1,v_2,\ldots,v_p\}$,  there are $\binom{r-p}{k-p}$ distinct $k$-tuples $w_1 < \cdots < w_k$ such that $\{v_1,\ldots,v_p\} \subset \{w_1,\ldots, w_k\} \subset \brr{\s}$.
Thus, summing first over $\rb\in \Lambda_{k,*}$ we have
$$C_{v_1, \ldots, v_p}(\bs) = \sum_{k = p}^r \binom{r-p}{k-p}(-1)^{k-p}c^{p-r}.$$
By the binomial theorem, $\sum_{k=p}^r (-1)^{k-p}\binom{r-p}{k-p}= (1-1)^{r-p}=0$. Thus $C_{v_1\cdots v_p}(\bs)=0$. \rev{Hence the only terms that contribute are indicators involving $p = r$ variables. 
The indicator $\ind\{X_{v_1}=\cdots =X_{v_r}\}$ shows up in the double sum only when $\rb=\brr{\mathbf{s}}$. This implies $p=r=k$, and $\{v_1,\ldots,v_r\}=\brr{\mathbf{s}}$. In that case, this forces the coefficient $C_{v_1,\ldots,v_r}(\mathbf{s})$ to $(-1)^{k-p}c^{p-r}=1$. Thus, in conclusion, we have 
\begin{align*}
    \sum_{k=2}^r\sum_{\rb\in \Lambda_{k,*}} \widetilde{Y}_{\rb}= \frac1{\Aut(H)}\sum_{\s\in V(G_n)_r} a_{H,\s} \ind\{X_{s_1}=\cdots=X_{s_r}\},
\end{align*}
which is precisely the definition of $T(H,G_n)$ from~\eqref{d:thgn}. Thus~\eqref{ob1} holds.}
\end{proof}
\section{Proof of Theorem~\ref{t:clt}}\label{sec3}

In this section, we prove Theorem~\ref{t:clt}, showing a quantitative central limit theorem for $Z(H, G_n)$.
We start by leveraging the decomposition  of $Z(H,G_n)$ from Proposition~\ref{p:zrgn0}, which allows us to write $Z(H,G_n)$ as sum of a martingale difference sequence $\{U_t\}_{t\in V(G_n)}$ (defined in~\eqref{d:ut}).
We then apply a following normal approximation result for martingale difference sequences of~\cite[Theorem, Eq (1)]{HEY70}, to find
\begin{align}\label{e:zrbd}
    \sup_{x\in \R}|\mathbb{P}[Z(H,G_{n})\leq x]-\Phi(x)| &\lesssim (A+B)^{1/5},
\end{align}
where
\begin{align}\label{d:ab}
A:=\sum_{t=1}^{v(G_n)}\mathbb{E}[U_{t}^{4}],\quad  B:=\Var\left[\sum_{t=1}^{v(G_n)}U_{t}^{2}\right].
\end{align}
This reduces our original problem to the task of understanding the quantities $A$ and $B$, which we do through careful subgraph counting.

\begin{lemma}\label{l:aub}
We have that 
\begin{align}\label{e:a0}
A \lesssim \frac{N(\gj, G_n)}{N(\mfh,G_n)^{2}} + \frac{N(\mth, G_n)^{1/4}}{N(\mfh,G_n)^{1/2}},
\end{align}
where $\gj$ and  $\mathcal{Q}_k^H$ denote the families of good joins and $2$-shared $k$-joins of $H$ defined in  Definitions~\ref{def:good} and~\ref{2shared} respectively.
\end{lemma}

\begin{lemma}\label{l:bub}
We have that 
\begin{align}\label{e:b0}
B \lesssim \frac{N(\gj, G_n)}{N(\mfh,G_n)^{2}} + \frac{N(\mth, G_n)^{1/4}}{N(\mfh,G_n)^{1/2}},
\end{align}
where $\gj$ and $\mathcal{Q}_k^H$ are good joins and $2$-shared $k$-joins of $H$ defined in  Definition~\ref{def:good} and Definition \ref{2shared} respectively.
\end{lemma}

We defer the proofs of the above pair of lemmas to Section~\ref{sec3.2} and Section~\ref{sec3.3} respectively. Assuming these lemmas, we can complete the proof of Theorem \ref{t:clt}.

\begin{proof}[Proof of Theorem \ref{t:clt}] The error bound in Theorem \ref{t:clt}, \eqref{e:clt}, follows immediately when $N(\gj,G_n) \ge N(\mfh,G_n)^2$. Assume instead that $N(\gj,G_n)<N(\mfh,G_n)^2$. Then 
\begin{align*}
    \frac{N(\gj, G_n)}{N(\mfh,G_n)^{2}} \le \frac{N(\gj, G_n)^{1/4}}{N(\mfh,G_n)^{1/2}}.
\end{align*}
Since $\mth \subset \gj$, we conclude $ N(\mth,G_n)\le N(\gj,G_n)$. This yields 
\begin{align*}
    \frac{N(\mth, G_n)^{1/4}}{N(\mfh,G_n)^{1/2}} \lesssim \frac{N(\gj, G_n)^{1/4}}{N(\mfh,G_n)^{1/2}}.
\end{align*}
Plugging these bounds into~\eqref{e:a0} and~\eqref{e:b0}, in view of the martingale central limit theorem bound~\eqref{e:zrbd}, we arrive at~\eqref{e:clt}. This completes the proof.
\end{proof}

It remains to prove Lemmas~\ref{l:aub} and~\ref{l:bub}.

\subsection{Preliminary observations} 
We begin with some preliminary observations about $Y_{\rb}$ (defined in~\eqref{e:defY}) which will prove useful in our subsequent analysis.

\begin{lemma}\label{l:yprop}
Fix $q\in \Z_{>0}$. Take $k_1,k_2,\ldots, k_q \in \ll 2,r\rr$ and recall $Y_{\rb}$  from \eqref{e:defY}. 
\begin{enumerate}[(a)]
    \item Let $\rb_j\in \Lambda_{k_j,*}$ for $j\in [q]$. Then,
    $\EE[\prod_{j=1}^q Y_{\rb_j}]$ is nonzero only if
    \begin{align}\label{e:subcond}
        \rb_i \subset \bigcup_{j\neq i} \rb_j \mbox{ for all } i\in [q].
    \end{align}
    
    \item Let $\rb_j\in \Lambda_{k_j,*}$ for $j\in [q]$.  Given $\alpha_i \in \Z_{>0}$, we have
    \begin{align*}
        \Ex\left[\prod_{j=1}^q Y_{\rb_j}^{\alpha_j}\right] \lesssim \prod_{j=1}^q D_{\rb_j}^{\alpha_j},
    \end{align*}
    where $D_{\rb}$ is defined in \eqref{d:dw}.
    Here the constant in $\lesssim$ may depend on $\alpha_j$'s, $q$ and $k$, but it does not depend $G_n$ and on our choices $\rb_1, \ldots, \rb_q$.

    \item Fix any $k,k'\in \ll 2,r\rr$. For any $\rb\in \Lambda_{k,*}, \rb' \in \Lambda_{k',*}$ and $v\in V(G_n)$, we have $\Ex[Y_{\rb}Y_{\rb'}\mid X_v]=\Ex[Y_{\rb}Y_{\rb'}]$.
\end{enumerate}
\end{lemma}
\begin{proof} \textbf{(a).} Suppose~\eqref{e:subcond} fails. Without loss of generality we may assume $\rb_1 \not\subset \cup_{j=2}^q \rb_j$. Thus there exists a vertex $v$ which appears in $\rb_1$ and does not appear in any $\rb_j$ for $j\in \ll 2,q\rr$. By the tower property of the expectation we have
    \begin{align*}
    \EE\left[\prod_{j=1}^q Y_{\rb_j}\right] &= \EE\left[\EE\left[\prod_{j=1}^q Y_{\rb_j} \mid \{X_i\}_{i \in [v(G_n)] \backslash \{v\}}\right]\right] =  \EE\left[\prod_{j=2}^q Y_{\rb_j}\EE\left[Y_{\rb_1} \mid \{X_i\}_{i \in [v(G_n)] \backslash \{v\}}\right]\right].  
    \end{align*}
    By Lemma~\ref{l:ycond}, the inner expectation in the right hand side of the above equation is zero, forcing $\Ex[\prod_{j=1}^q Y_{\rb_j}]$ to be zero as well. This proves part (a).

\medskip    

\noindent \textbf{(b).} For $\rb\in \Lambda_{k,*}$, we denote 
$${M}_{\rb}(p):=\sum_{1\le i_1<\cdots<i_p\le k}\left[\frac1{c^{r-p}}\ind\{X_{w_{i_1}}=\cdots=X_{w_{i_p}}\}-\frac{1}{c^{r-1}}\right],$$
so that ${Y}_{\rb}:=D_{\rb}\sum_{p=2}^{k}(-1)^{k-p}M_{\rb}(p)$. In particular,
\rev{\begin{align}\label{e:calct}
   \Ex\left[\prod_{j=1}^q Y_{\rb_j}^{\alpha_j}\right] \le \left[\prod_{j=1}^q D_{\rb_j}^{\alpha_j} \right]\cdot \Ex\left[\prod_{j=1}^q\left|\sum_{p=2}^{k_j} (-1)^{k_j-p}M_{\rb_j}(p)\right|^{\alpha_j}\right],
\end{align}
For the expectation term above, applying H\"older inequality we get
\begin{align*}
    \Ex\left[\prod_{j=1}^q\left|\sum_{p=2}^{k_j} (-1)^{k_j-p}M_{\rb_j}(p)\right|^{\alpha_j}\right] & \le \prod_{j=1}^q \Ex\left[\left|\sum_{p=2}^{k_j} (-1)^{k_j-p}M_{\rb_j}(p)\right|^{q\alpha_j}\right]^{\frac1q}   \lesssim \prod_{j=1}^q \left[ \sum_{p=2}^{k_j} \Ex\left[\left|M_{\rb_j}(p)\right|^{q\alpha_j}\right]\right]^{\frac1q}.
\end{align*}
Here  the constant in $\lesssim$ depends only $\alpha_j$'s, $c$, $q$ and $k$. The last inequality follows from the standard inequality $|\sum_{p=2}^k r_p|^m \lesssim_{k,m} \sum_{p=2}^k |r_p|^m$ for $m\ge 1$. For the inner expectation above, following the definition of $M_{\rb}(p)$, we may again apply this standard inequality to get 
\begin{align*}
    \Ex\left[\left|M_{\rb}(p)\right|^{q\alpha}\right] \lesssim \sum_{1\le i_1<\cdots<i_p\le k} \Ex \left[{c^{p-r}}\ind\{X_{w_{i_1}}=\cdots =X_{w_{i_p}}\}\right]^{q\alpha} \le k^p c^{(p-r)q\alpha}. 
\end{align*}
The above expectation estimate leads to
\begin{align*}
    \Ex\left[\prod_{j=1}^q\left|\sum_{p=2}^{k_j} (-1)^{k_j-p}M_{\rb_j}(p)\right|^{\alpha_j}\right] \lesssim 1
\end{align*}
where the constant in $\lesssim$ depends only $\alpha_j$'s, $c$, $q$ and $k$. Inserting this back in \eqref{e:calct}} leads to the desired estimate in part (b). 

\medskip

\noindent \textbf{(c).}  
Using the same notation as above, observe that 
\begin{align*}
\EE[Y_{\rb} Y_{\rb'} \mid X_v] &= D_{\rb} D_{\rb'} \sum_{p=2}^{k} \sum_{q=2}^{k} (-1)^{2k-p-q} \EE\left[ M_{\rb}(p)M_{\rb'}(q) \mid X_v \right].
\end{align*}
A product of indicators is still an indicator, so $M_{\rb}(p)M_{\rb}(q)$ is just a linear combination of indicators of the form $\ind\{{X_{v_1}=\cdots=X_{v_p}}\}$ where $\{v_1,v_2,\ldots,v_p\}\subset \brr{\rb}\cup \overline{\rb}'$. It suffices to check
\begin{align}\label{e:check}
    \P[X_{v_1}=X_{v_2}=\cdots=X_{v_p}\mid X_v]=\P[X_{v_1}=X_{v_2}=\cdots=X_{v_p}].
\end{align}
If $v\not\in \{v_1,\ldots,v_p\}$, the above equality is immediate. If $v\in \{v_1,\ldots,v_p\}$. The left hand side requires all the other vertices to have same color as $X_v$. This happens with probability $c^{1-p}$. We can compute the right hand side:
\begin{align*}
    \P[X_{v_1}=X_{v_2}=\cdots=X_{v_p}]=\sum_{a=1}^c \P[X_{v_1}=X_{v_2}=\cdots=X_{v_p}=a]=\sum_{a=1}^c c^{-p}=c^{1-p}.
\end{align*}
This proves \eqref{e:check}, proving (c).
\end{proof}

We also give an estimate for the variance $\sigma_H^2$ defined in \eqref{d:zhgn} in terms of graph counts.

\begin{lemma}\label{sigbd} For any fixed $c \ge 2$, 
\begin{align}
    \label{e:sigbd}
    \sigma_H^2  \asymp N(\mfh, G_n)\asymp \sum_{k=2}^r \sum_{\rb\in \Lambda_{k,*}} D_{\rb}^2,
\end{align}
 where $D_{\rb}$ is defined in \eqref{d:dw} and $\mfh$ is defined in Definition \ref{2shared}.
\end{lemma}
\begin{proof} From \eqref{d:thgn} we have
\begin{align*}
    \sigma_H^2= \vr[T(H,G_n)] & = \frac1{|\Aut(H)|^2}\sum_{\s_1,\s_2\in V(G_n)_r} a_{H,\s_1}a_{H,\s_2}\Cov[\ind\{{X}_{=\bs_1}\},\ind\{{X}_{=\bs_2}\}].
\end{align*}
 Clearly the covariance $\Cov[\ind\{{X}_{=\bs_1}\},\ind\{{X}_{=\bs_2}\}]$ is zero when $|\brr\bs_1\cap\brr\s_2|=0$ due to independence. If $|\brr\bs_1\cap\brr\s_2|\ge 1$, then 
$$\Cov[\ind\{{X}_{=\bs_1}\},\ind\{{X}_{=\bs_2}\}]=\Ex[\ind\{X_{=\bs_1}\}\ind\{X_{=\bs_2}\}]-\P[X_{=\bs_1}]\P[X_{=\bs_2}]=c^{1-|\brr\s_1\cup \brr\s_2|}-c^{2-|\brr\s_1|-|\brr\s_2|}.$$ 
Thus the covariance is still zero when $|\brr\s_1\cap\brr\s_2|=1$. For $|\brr\s_1\cap \brr\s_2|\ge 2$ 
observe that
$$ c^{-2r}(1-c) \le \Cov[\ind\{{X}_{=\bs_1}\},\ind\{{X}_{=\bs_2}\}]=c^{1-|\brr\s_1\cup \brr\s_2|}-c^{2-|\brr\s_1|-|\brr\s_2|} \le c^{1-r}.$$ 
Thus,
$$\sigma_H^2 \asymp \sum_{\substack{\s_1,\s_2\in V(G_n)_r \\ |\brr\s_1\cup\brr\s_2|\ge 2}} a_{H,\bs_1}a_{H,\bs_2} \asymp N(\mfh,G_n).$$
This proves the first part in \eqref{e:sigbd}. For the second part, for all $k\in \ll2,r\rr$, and $\rb\in \Lambda_{k,*}$, we define
$$R_{\rb}:= \sqrt{\sum_{\substack{\s_1,\s_2\in V(G_n)_r \\ \brr\s_1\cup\brr\s_2=\rb}} a_{H,\bs_1}a_{H,\bs_2}},$$ 
Clearly $R_{\rb}^2\le D_{\rb}^2$. We thus have
\begin{align*}
    \sigma_H^2 \asymp \sum_{\substack{\s_1,\s_2\in V(G_n)_r \\ |\brr\s_1\cup\brr\s_2|\ge 2}} a_{H,\bs_1}a_{H,\bs_2} \asymp \sum_{k=2}^r\sum_{\rb\in \Lambda_{k,*}} R_{\rb}^2 \le \sum_{k=2}^r\sum_{\rb\in \Lambda_{k,*}} D_{\rb}^2.
\end{align*} 
To show the other side inequality observe that
\begin{align*}
    \sum_{\rb\in \Lambda_{k,*}} D_{\rb}^2 = \sum_{\rb\in \Lambda_{k,*}} \sum_{\ell=k}^r\sum_{\rb'\in \Lambda_{\ell,*}, \rb'\supset \rb} R_{\rb'}^2 = \sum_{\ell=k}^r\sum_{\rb'\in \Lambda_{\ell,*}} \binom{\ell}{k}R_{\rb'}^2 \le 2^r\sum_{\ell=k}^r\sum_{\rb'\in \Lambda_{\ell,*}} R_{\rb'}^2.
\end{align*}
Thus,
\begin{align*}
    \sum_{k=2}^r\sum_{\rb\in \Lambda_{k,*}} D_{\rb}^2 \lesssim \sum_{k=2}^r\sum_{\rb\in \Lambda_{k,*}} R_{\rb}^2 \asymp \sigma_H^2.
\end{align*}
This proves the other side inequality, completing the proof of Lemma \ref{sigbd}.
\end{proof}

While upper bounding the terms $A$ and $B$ defined in \eqref{d:ab}, we encounter error functions involving the following expression
\begin{align}\label{d:skgn}
    \skr(G_n) := \sum_{2\le k_1,k_2\le r}\sum_{t=1}^{v(G_n)} \sum_{(\rb_1,\rb_2)\in \Lambda_{k_1,k_2,t}^{(1)}} D_{\rb_1}^2 D_{\rb_2}^2,
\end{align}
where $\Lambda_{k_1,k_2,t}^{(1)}$ and $D_{\rb}$ are defined in \eqref{e:ovlap} and \eqref{d:dw} respectively. Observe that $\skr(G_n)$ as defined above has a dependence on our choice of vertex ordering. Through the rest of the article, we will fix a particular ordering on the vertices of $G_n$, per the subsequent Lemma~\ref{l:skgn}.

Before going on to prove Lemma~\ref{l:aub} and Lemma~\ref{l:bub}, we end this subsection with a technical lemma that provides an upper bound for this error function $\skr(G_n)$.

\begin{lemma}\label{l:skgn} 
There exists an ordering of the vertices of $V(G_n)$ such that for $\skr(G_n)$ defined as in~\eqref{d:skgn} with respect to this ordering,
\begin{align}\label{e:ubd}
\skr(G_n) \lesssim N(\mfh,G_n)^{3/2} \cdot N(\mth, G_n)^{1/4},
\end{align}
where $\mathcal{Q}_k^H$ is the collection of all $2$-shared $k$-joins defined in Definition \ref{2shared}.
\end{lemma}

\begin{proof} 
Recall the definition of a strictly increasing tuple and the associated Definition~\ref{d:skt}.
 We claim that 
\begin{align}
    \label{ubd}
    \skr(G_n) \lesssim \sum_{1\le a<b<t\le v(G_n)} D_{at}^2D_{bt}^2,
\end{align}
where $D_{at} = D_{\bw}$ for $\bw = (a, t)$.
To see this, we start by fixing $k\in \ll 2,r\rr$ and $t\in V(G_n)$. 
We make the following initial observation. By summing over $D_{\rb}$ where $\rb = (w_1, \ldots, w_{k-2}, y, t) \in \Lambda_{k, t}$ for some fixed $y < t$, we find that 
\begin{align}\label{sk1}\sum_{\rb \in \Lambda_{k, t} : \rb_{k-1} = y} D_{\rb}^2 \le \left(\sum_{\rb \in \Lambda_{k, t} : \rb_{k-1} = y} D_{\rb} \right)^2 \lesssim D_{ y t }^2, 
\end{align}
where the first inequality uses the fact that $D_{\rb} \ge 0$ and the second follows from the definition of $D_{\rb}$ as counting the number of copies of $H$ that pass through $\rb$.

Recall that $(\rb_1,\rb_2)\in \Lambda_{k_1,k_2,t}^{(1)}$, then $\rb_1\in \Lambda_{k_1,t}$, $\rb_2\in \Lambda_{k_2,t}$, and $\rb_1\cap \rb_2 =\{t\}.$ We write $\rb_{1}=(w_{1,1},\ldots,w_{1,k_{1}})$ and $\rb_{2}=(w_{2,1},\ldots,w_{2,k_{2}})$. Consequently, each term of the innermost sum in $\skr(G_n)$ can be bounded as follows:
\begin{align*}
    \sum_{(\rb_1,\rb_2)\in \Lambda_{k_1,k_2,t}^{(1)}} D_{\rb_1}^2 D_{\rb_2}^2 \lesssim  \sum_{\substack{\rb_i,\in \Lambda_{k_i,t}, i=1,2 \\ w_{1,k_1-1}\neq w_{2,k_2-1}}} D_{\rb_1}^2 D_{\rb_2}^2 & \lesssim \sum_{\substack{w_{1,k_1-1}, w_{2,k_2-1} \\ w_{1,k_1-1}\neq w_{2,k_1-1}, w_{i,k_i-1}<t}}  D_{w_{1,k_1-1},t}^2D_{w_{2,k_2-1},t}^2. 
\end{align*}
Above, the first inequality follows by extending the sum over all pairs of tuples $\rb_1,\rb_2$ which are required to have distinct penultimate entries $w_{1,k_{1}-1}\neq w_{2,k_{2}-1}$. The next bound follows by applying~\eqref{sk1}; we sum over the first $k_1-2$ and $k_2-2$ entries of $\rb_1$ and $\rb_2$ respectively. The right hand side of the above inequality is precisely equal to the right hand side of \eqref{ubd}. This proves~\eqref{ubd}. 

\medskip

The rest of the proof is adapted from Lemma 3.1 in \cite{BHA22} which deals with the case when $H$ is a triangle. We include this argument in order to be self-contained.
Following ~\cite{BHA22}, we use the following notation for the rest of the argument.
\begin{itemize} 
    \item $\mathcal{Q}^{H}_{2,s}$ is the collection of $2$-joins of $H$ in which the copies share at least two vertices with $s$ as one of the vertices
    \item $\mathcal{Q}^{H}_{2, s_{1},s_{2}}$ is the collection of $2$-joins of $H$ in which the copies share vertices $s_{1}$ and $s_{2}$
    \item $\mathcal{L}_{s_{1},s_{2},s_{3}}$ is the collection of $4$-joins of $H$ formed by two copies which share vertices $s_{1}$ and $s_{3}$ and two copies which share the vertices $s_{2}$ and $s_{3}$.
\end{itemize}
We set $\xi_{t}:=N(\mathcal{Q}_{2,t}^{H},G_{n})$ and fix an ordering of the vertices so that \begin{align}\label{ord}
    \xi_{1}\geq\cdot\cdot\cdot\geq\xi_{v(G_{n})}.
\end{align} The expression $D_{at}^{2}D_{bt}^{2}$ counts the possible combinations of two copies of $H$ with $\{a,t\}$ as shared vertices, and two copies of $H$ with $\{b,t\}$ as an shared vertices. Clearly these four copies of $H$ together form one of the motifs from the collection $\mathcal{L}_{a,b,t}.$ 
Thus, 
\begin{align}\label{ubd1}
    \sum_{1\le a<b<t\le v(G_n)} D_{at}^2D_{bt}^2 & \lesssim \sum_{1\leq a<b<t\leq v(G_{n})} N(\mathcal{L}_{a,b,t},G_{n}).
\end{align} 
Fix a vertex $a \ge 1$. If we fix $t > a$, an element of $\mcL_{a,b,t}$ in $G_n$ (for some $b$ satisfying $a < b < t$) arises from a choice of two copies of $H$ that share vertices $a$ and $t$ (there are at most $N(\mathcal{Q}_{2,a,t}^H, G_n)$ such choices) and a pair of copies of $H$ sharing vertices $b, t$ (at most $\xi_t$ choices over all possibilities of $b$). This implies the following upper bound:
\begin{align*}
    \sum_{b,t: a<b<t\leq v(G_{n})}N(\mathcal{L}_{a,b,t},G_{n}) \lesssim \sum_{t: a<t\leq v(G_{n})} N(\mathcal{Q}^{H}_{2,a,t},G_{n})\cdot \xi_{t}.
\end{align*} 

By a similar analysis, we can also obtain the following different upper bound:
\begin{align*}
    \sum_{b,t:1\leq a<b<t\leq v(G_{n})}N(\mathcal{L}_{a,b,t},G_{n}) & \lesssim \xi_{t}\xi_{a}\leq \xi_{a}^{2}.
\end{align*} 
The above inequalities follow by noting that we can choose two copies of $H$ sharing two vertices, one of which is $a$, in $\xi_{a}$ ways, and two copies of $H$ sharing two vertices, one of which is $t$, in $\xi_t\le \xi_a$ ways. By taking a geometric average of the two bounds and summing over $a$ we find that 
\begin{align}
    \sum_{1\leq a<b<t\leq v(G_{n})}N(\mathcal{L}_{a,b,t},G_{n}) \lesssim \sum_{a=1}^{v(G_{n})}\left(\sqrt{\xi_{a}^{2}}\left(\sum_{t=1}^{v(G_{n})}N(\mathcal{Q}_{2,a,t}^{H},G_{n})\xi_{t}\right)^{1/2}\right). \label{upl}
\end{align}
Observe that $N(\mathcal{Q}_4^H,G_n) \gtrsim \max_{a,t} N(\mathcal{Q}_{2,a,t}^H)^2$ (recall that $\mathcal{Q}_k^H$ is the collection of all $2$-shared $k$-joins of $H$, $k$-joins of $H$ with $2$ vertices common to all copies). We also have $\sum_{t} \xi_t \lesssim N(\mfh,G_n)$. Thus,
\begin{align*}
    \mbox{r.h.s.~of \eqref{upl}} \lesssim N(\mth,G_n)^{1/4}N(\mfh,G_n)^{3/2}.
\end{align*}
In view of \eqref{ubd} and \eqref{ubd1}, we get the desired bound in \eqref{e:ubd}. This completes the proof.
\end{proof}

In the rest of this section, we fix an ordering on $V(G_n)$ given by \eqref{ord} in Lemma~\ref{l:skgn}.

\subsection{Proof of Lemma~\ref{l:aub}}\label{sec3.2} 
Recall $U_t$ from \eqref{d:ut}.
Using the fact
$(\sum_{k=2}^p a_k)^4 \lesssim_p \sum_{k=2}^p a_k^4$, we have
\begin{align}
    \label{e:a1}
    A := \sum_{t=1}^{v(G_n)}\mathbb{E}[U_{t}^{4}] \lesssim \frac{1}{\sigma_H^4}  \sum_{t = 1}^{v(G_n)} \sum_{k = 2}^r \EE\left[\left(\sum_{\rb\in \Lambda_{k,t}} Y_{\rb} \right)^4 \right],
\end{align}
where $Y_{\rb}$ is defined in~\eqref{e:defY}. We now expand the fourth powers on the right hand side of \eqref{e:a1} to get 
\begin{align}\label{e:a2}
\mbox{r.h.s.~of \eqref{e:a1}}
\lesssim \frac1{\sigma_H^4}\sum_{k=2}^r (A_1^{(k)} + A_2^{(k)} + A_3^{(k)} + A_4^{(k)}),
\end{align}
where
\begin{align}\label{def:a12}
A_1^{(k)} &:=  \sum_{t = 1}^{v(G_n)} \sum_{\rb\in \Lambda_{k,t}} \EE\left[Y_{\rb}^4\right], \quad
A_2^{(k)} :=  \sum_{t = 1}^{v(G_n)} \sum_{\substack{\rb_1, \rb_2 \in \Lambda_{k,t} \\ \rb_1\neq \rb_2}}\EE\left[Y_{\rb_1}^2 Y_{\rb_2}^2\right],
\end{align}
\begin{align}\label{def:a3}
A_3^{(k)} &:=  \sum_{t = 1}^{v(G_n)}   \sum_{\substack{\rb_1,\rb_2,\rb_3 \in \Lambda_{k,t} \\ \rb_i \neq \rb_j \ \operatorname{for} \ i\neq j}}\EE\left[Y_{\rb_1}^2 Y_{\rb_2} Y_{\rb_3}\right],
\end{align}
\begin{align}\label{def:a4}
A_4^{(k)} &:=  \sum_{t = 1}^{v(G_n)} \sum_{\substack{\rb_1,\rb_2,\rb_3,\rb_4 \in \Lambda_{k,t} \\ \rb_i \neq \rb_j \ \operatorname{for} \ i\neq j}}\EE\left[Y_{\rb_1} Y_{\rb_2} Y_{\rb_3} Y_{\rb_4}\right].
\end{align}
Observe that terms of the form $Y_{\rb_1}^3Y_{\rb_2}$ do not appear on the right hand side of \eqref{e:a2} as whenever $\rb_1\neq \rb_2$, we have $\Ex[Y_{\rb_1}^3Y_{\rb_2}\mid X_v, v\in \rb_1]=Y_{\rb_1}^3\Ex[Y_{\rb_2}\mid X_v, v\in \rb_1]=0$, by Lemma~\ref{l:ycond}.

The following claim provides an estimate for each $A_i^{(k)}$.

\begin{claim}\label{c:aest} For each $k\in \ll 2,r\rr$ we have
\begin{align}\label{e:cl0}
A_2^{(k)}  \lesssim \sum_{t=1}^{v(G_n)}\sum_{\rb_1,\rb_2\in \Lambda_{k,t}} D_{\rb_1}^2D_{\rb_2}^2 \lesssim \skr(G_n) +  N(\gj, G_n),
\end{align}
and
\begin{align*}
\quad A_1^{(k)}, A_3^{(k)},
    A_4^{(k)} \lesssim  N(\gj, G_n),
\end{align*}
where $\skr(G_n)$ and $\gj$ are defined in Definition~\ref{d:skgn} and Definition~\ref{def:good} respectively.  
\end{claim}
Assuming the claim, we can bound the right hand side of \eqref{e:a2} as follows
\begin{align*}
\mbox{r.h.s.~of \eqref{e:a2}}
\lesssim \frac1{\sigma_H^4}  \left(N(\gj, G_n) + \skr(G_n) \right). 
\end{align*}
Plugging in the estimate for $\sigma_H^2$ from Lemma~\ref{sigbd} and the upper bound for $\skr(G_n)$ from Lemma~\ref{l:skgn} into the right hand side of the equation above, we arrive at \eqref{e:a0}. This completes the proof of Lemma~\ref{l:aub} except for the proof of Claim~\ref{c:aest}.

\begin{proof}[Proof of Claim~\ref{c:aest}] \underline{\textbf{$A_1^{(k)}$ terms.}} Recall $A_1^{(k)}$ from \eqref{def:a12} and $\Lambda_{k,*}$, $\Lambda_{k,t}$ from Definition~\ref{d:skt}. We see that
\begin{align*}
A_1^{(k)}=\sum_{t=1}^{v(G_n)}\sum_{\rb\in\Lambda_{k,t}} \EE\left[ {Y}_{\rb}^4 \right]= \sum_{\rb\in\Lambda_{k,*}} \EE\left[ {Y}_{\rb}^4 \right] 
\overset{\text{Lem.~\ref{l:yprop}(b)}}\lesssim \sum_{\rb\in \Lambda_{k,*}} D_{\rb}^4 
\lesssim N(\gj, G_n).
\end{align*}
The last inequality follows from the observation that each contribution to $D_{\rb}^4$ is given by a quadruple of copies $H$ that have at least $k$ vertices amongst all four copies. For $k\ge 2$, such a join belongs to $\gj$ (see Definition~\ref{def:good} and~\ref{2shared}).

\medskip

\noindent\underline{\textbf{$A_2^{(k)}$ terms.}} Recall $A_2^{(k)}$ from \eqref{def:a12}. Applying Lemma~\ref{l:yprop}(b) again, we find that
\begin{align*}
A_2^{(k)}  :=\sum_{t=1}^{v(G_n)} \sum_{\rb_1,\rb_2\in \Lambda_{k,t}}\EE\left[Y_{\rb_1}^2 Y_{\rb_2}^2\right] & \lesssim \sum_{t=1}^{v(G_n)} \sum_{\rb_1,\rb_2\in \Lambda_{k,t}}D_{\rb_1}^2 D_{\rb_2}^2 \\ & = \sum_{t=1}^{v(G_n)}\sum_{\substack{\rb_1,\rb_2\in \Lambda_{k,t} \\ {\rb}_1 \cap {\rb}_2 = \{t\}}}D_{\rb_1}^2 D_{\rb_2}^2 +\sum_{t=1}^{v(G_n)} \sum_{\substack{\rb_1,\rb_2\in \Lambda_{k,t} \\ |{\rb}_1 \cap {\rb}_2|\ge 2}}D_{\rb_1}^2 D_{\rb_2}^2.
\end{align*}
The first term in the resulting equation is a subset of terms contributing to $\skr(G_n)$ (see~\eqref{d:skgn}). The second term counts the number of $4$-joins of $H$, where there are at least two vertices common to all 4 copies. Such a $4$-join falls into $\gj$. This proves~\eqref{e:cl0}.

\medskip

\noindent\underline{\textbf{$A_3^{(k)}$ terms.}} For $A_3^{(k)}$ defined in \eqref{def:a3}, consider any individual term $\EE\left[Y_{\rb_1}^2 Y_{\rb_2} Y_{\rb_3}\right]$, where $\rb_1,\rb_2,\rb_3 \in \Lambda_{k,t}$. By the definition of $\Lambda_{k,t}$ (see~\eqref{e:lkt}), $\rb_1,\rb_2,\rb_3$ must all share common vertex $t$. Thus, $\left|{\rb}_1 \cap \left\{{\rb}_2 \cup {\rb}_3 \right\} \right| \ge 1.$
We claim that $\EE\left[Y_{\rb_1}^2 Y_{\rb_2} Y_{\rb_3}\right]$ is nonzero only if $\left|{\rb}_1 \cap \left\{{\rb}_2 \cup {\rb}_3 \right\} \right| \ge 2.$ 
Indeed if ${\rb}_1 \cap \left\{{\rb}_2 \cup {\rb}_3 \right\}=\{t\}$, then   $$\EE\left[Y_{\rb_1}^2 Y_{\rb_2} Y_{\rb_3}\mid X_{v}, v\in \rb_1 \right] = Y_{\rb_1}^2\EE\left[ Y_{\rb_2}Y_{\rb_3}\mid X_t\right].$$ 
\rev{As $t\in \rb_3$, by tower property of the conditional expectation we have
\begin{align*}
    \EE\left[ Y_{\rb_2}Y_{\rb_3}\mid X_t\right]=\EE\left[ Y_{\rb_3} \EE\left[ Y_{\rb_2}\mid X_v, v\in \rb_3 \right] \mid X_t \right].
\end{align*}
But, also by Lemma~\ref{l:ycond}, $\EE\left[ Y_{\rb_2}\mid X_v, v\in \rb_3\right]=0$ as $\rb_2\neq \rb_3$.} This leads to a contradiction.  Appealing to Lemma~\ref{l:yprop}(b) we have
$$A_3^{(k)} \lesssim \sum_{\substack{\rb_1,\rb_2,\rb_3\in \Lambda_{k,*} \\ \left|{\rb}_1 \cap \left\{{\rb}_2 \cup {\rb}_3 \right\} \right| \ge 2}} D_{\rb_1}^2 D_{\rb_2} D_{\rb_3}.$$ The right hand side of the above equation counts number of $4$ joins of $H$ such that there are two vertices, say $v_1,v_2$, where $v_1$ appears in all $4$ copies of $H$ and $v_2$ appears in at least three of these copies. This forces conditions~\eqref{good} and~\eqref{vgood} to be satisfied. Thus the right hand side of the above equation is bounded  above by $N(\gj,G_n)$ (up to a constant), which is precisely what we wanted to show.  

\medskip

\noindent\underline{\textbf{$A_4^{(k)}$ terms.}} Finally, for $A_4^{(k)}$ (defined in \eqref{def:a4}), we again consider a representative term in the sum. By Lemma~\ref{l:yprop}(a), 
$\EE\left[Y_{\rb_1} Y_{\rb_2} Y_{\rb_3} Y_{\rb_4}\right]$ is nonzero only if \begin{align}
    \label{e:a3}
    {\rb}_i \subset \bigcup_{j \neq i} {\rb}_j \mbox{ for all } i \in [4].
\end{align}  Take any $4$-join of $H$, say $H_1\cup H_2 \cup H_3\cup H_4$, with $\rb_i\subset V(H_i)$. Any such $4$-join satisfies \eqref{good}. We claim that these $4$-joins are good joins, which means they must also satisfy \eqref{vgood}. If they do not satisfy \eqref{vgood}, then (permuting if necessary) 
$$v(H_1\cup H_2 \cup H_3\cup H_4)=v(H_1\cup H_2 )+v(H_3\cup H_4)-1.$$ 
Since $t$ is a common vertex for all $H_i$, we must have $\{{\rb}_1\cup {\rb}_2\} \cap \{{\rb}_3\cup {\rb}_4\} = \{t\}.$ 
By the above condition, there must be some vertex $v \in \rb_1 \backslash \{t\}$ which does not appear in ${\rb}_2\setminus \{t\}$. By~\eqref{e:a3}, such a vertex must appear in $\{{\rb}_3\cup {\rb}_4\}$. This is a contradiction.

Appealing to Lemma~\ref{l:yprop}(b) again and summing over all good joins we obtain the desired upper bound for $A_4^{(k)}$.
\end{proof}
\subsection{Proof of Lemma~\ref{l:bub}} \label{sec3.3}
Recall from \eqref{d:ab} that $B$ is the variance of the sum of $U_t^2$ for $U_t$ as defined in~\eqref{d:ut}. We begin by expanding $\sigma_{H}^{2}U_{t}^{2}$.
\begin{align*}
    \sigma_H^2 U_{t}^{2}  =\left(\sum_{k=2}^r \sum_{\rb\in \Lambda_{k,t}} Y_{\rb}\right)^2 & =\sum_{2\le k_1,k_2\le r} \sum_{\substack{\rb_1\in \Lambda_{k_1,t} \\ \rb_2\in \Lambda_{k_2,t}}} Y_{\rb_1}Y_{\rb_2} =\sum_{2\le k_1,k_2\le r} \sum_{u=1}^{k_1\wedge k_2} \sum_{\substack{(\rb_1,\rb_2)\in \Lambda_{k_1,k_2,t}^{(u)}}} Y_{\rb_1}Y_{\rb_2}.
\end{align*}
The notation $\Lambda_{k_1,k_2,t}^{(u)}$ was introduced in Definition~\ref{d:skt}. The last equality above follows by partitioning the inner sum into terms according to number of vertices $u$ on which $\rb_1$ and $\rb_2$ overlap. Since $\vr[\sum_{i=1}^p R_i] \lesssim_p \sum_{i=1}^p\vr[R_i]$, we find

\begin{align}\label{d:bk}
B & := \Var\left(\sum_{t=1}^{v(G_n)}U_{t}^{2}\right) \lesssim \frac1{\sigma_H^4}\sum_{2 \le k_1,k_2 \le r}\sum_{u=1}^{k_1\wedge k_2} B_{k_1,k_2}^{(u)}.
\end{align}
We define the term $B_{k_{1},k_{2}}^{(u)}$ as
\begin{equation}\label{d:bkk}
    \begin{aligned}
B_{k_1, k_2}^{(u)} & :=  \Var\bigg[\sum_{t=1}^{v(G_n)}  \sum_{(\rb_1,\rb_2)\in \Lambda_{k_1,k_2,t}^{(u)}} Y_{\rb_1}Y_{\rb_2} \bigg] \\ & =  \sum_{t_1,t_2\in V(G_n)}   \sum_{\substack{(\rb_1,\rb_2)\in \Lambda_{k_1,k_2,t_1}^{(u)} \\   (\rb_3,\rb_4)\in \Lambda_{k_1,k_2,t_2}^{(u)}}} \Cov[Y_{\rb_1}Y_{\rb_2},Y_{\rb_3}Y_{\rb_4}].
\end{aligned}
\end{equation}

We now proceed to estimate $B_{k_1,k_2}^{(u)}$ for different choices of $k_1,k_2,$ and $u$.

\begin{claim}\label{e:b1}  Recall $\skr(G_n)$ and $\gj$ from  Definitions~\ref{d:skgn} and~\ref{def:good} respectively. 
\begin{enumerate}[(a).]
    \item For each $k\in \ll 2,r\rr$ we have
    $B_{k,k}^{(k)} \lesssim N(\gj,G_n)$.
    \item Fix any $k_1,k_2\in \ll 2,r\rr$. For $u\in \ll 1,k_1\wedge k_2-1\rr$ or $k_1\neq k_2$ we have $$B_{k_1,k_2}^{(u)} \lesssim \skr(G_n)+N(\gj,G_n).$$
\end{enumerate}
\end{claim}

Plugging in the above estimates into~\eqref{d:bk}, we have
\begin{align*}
\mbox{r.h.s.~of \eqref{d:bk}}\lesssim \frac1{\sigma_H^4}  \left(N(\gj, G_n) + \skr(G_n) \right). 
\end{align*}
Finally, applying our earlier estimate for $\sigma_H^2$ from Lemma~\ref{sigbd} and the upper bound on $\skr(G_n)$ from Lemma~\ref{l:skgn} to the right hand side of above equation, we arrive at \eqref{e:a0}. This completes the proof of Lemma~\ref{l:bub} except for Claim~\ref{e:b1}.
\begin{proof}[Proof of Claim~\ref{e:b1}] \textbf{(a).}
We study the individual terms that appear in the characterization of $B_{k,k}^{(k)}$ of~\eqref{d:bkk}. For any $(\rb_1,\rb_2)\in \Lambda_{k,k,t_1}^{(k)}$ and  $(\rb_3,\rb_4)\in \Lambda_{k,k,t_2}^{(k)}$ we must have $\rb_1=\rb_2$ and $\rb_3=\rb_4$ (per~\eqref{e:ovlap}, since $u = k$). We upper bound the expression $\Cov[Y_{\rb_1}Y_{\rb_2},Y_{\rb_3}Y_{\rb_4}]=\Cov[Y_{\rb_1}^2,Y_{\rb_3}^2]$ in~\eqref{d:bkk} by considering three cases:
\begin{itemize}
    \item If $\rb_1,\rb_3$ are disjoint, then the covariance is zero by independence.
    \item If $\rb_1,\rb_3$ have precisely one vertex in common, say $v$ then by Lemma~\ref{l:yprop}(c), we have $$\Ex\left[Y_{\rb_1}^2Y_{\rb_3}^2\mid X_v\right]=\Ex\left[Y_{\rb_1}^2\mid X_v\right]\Ex\left[Y_{\rb_3}^2\mid X_v\right]=\Ex\left[Y_{\rb_1}^2\right]\Ex\left[Y_{\rb_3}^2\right],$$
    which forces the covariance of $Y_{\rb_1}^2$ and $Y_{\rb_3}^2$ to be zero.
    
    \item If $\rb_1,\rb_3$ have at least two vertices in common, then by Lemma~\ref{l:yprop} (b) we have $$\Cov[Y_{\rb_1}^2,Y_{\rb_3}^2]\le \Ex[Y_{\rb_1}^2Y_{\rb_3}^2] \lesssim D_{\rb_1}^2D_{\rb_3}^2.$$
    where $D_{\rb}$ is defined in \eqref{d:dw}.
\end{itemize}
Only the final list item will have nonzero contribution to the covariance bound. The expression $D_{\rb_1}^2D_{\rb_3}^2$ counts the number of $4$-joins of $H$ such that at least two copies of $H$ pass through $\rb_1$ and at least two copies of $H$ through $\rb_3$. Since $\rb_1$ and $\rb_3$ share at least two vertices, the corresponding $4$-joins satisfy~\eqref{good} and~\eqref{vgood} and thus belong to the collection $\gj$, defined in Definition~\ref{def:good}. By summing first over all such $\rb_1$ and $\rb_3$ and then over $t_1,t_2$, we arrive at Claim~\ref{e:b1}(a).

\medskip

\noindent\textbf{(b)}. Recall $\Lambda_{k_1,k_2,t}^{(u)}$ from \eqref{e:ovlap}. Take any $(\rb_1,\rb_2),\in \Lambda_{k_1,k_2,t_1}^{(u)}, (\rb_3,\rb_4)\in \Lambda_{k_1,k_2,t_2}^{(u)}$ meeting the conditions of part (b). Since $u<k_1\wedge k_2$ or $k_1\neq k_2$, there must be at least one vertex which does not appear in both $\rb_1$ and $\rb_2$. This forces $\Ex[Y_{\rb_1}Y_{\rb_2}]=0$ by Lemma~\ref{l:yprop}(a). Thus, $$\Cov[Y_{\rb_1}Y_{\rb_2},Y_{\rb_3}Y_{\rb_4}]=\Ex[Y_{\rb_1}Y_{\rb_2}Y_{\rb_3}Y_{\rb_4}].$$
By Lemma~\ref{l:yprop}(a), 
$\EE\left[Y_{\rb_1} Y_{\rb_2} Y_{\rb_3} Y_{\rb_4}\right]$ is nonzero only if \begin{align}
    \label{e:b}
    {\rb}_i \subset \bigcup_{j \neq i} {\rb}_j \mbox{ for all } i \in [4].
\end{align}
We upper bound $\EE\left[Y_{\rb_1} Y_{\rb_2} Y_{\rb_3} Y_{\rb_4}\right]$ in cases.

\smallskip

 \noindent{\textbf{Case 1.}} Suppose $t_1=t_2=t$. Let
 \begin{equation}\label{adef}
     \begin{aligned}
     \mathcal{A} & :=\left\{(\rb_1,\rb_2,\rb_3,\rb_4) \ \big| \ (\rb_1,\rb_2), (\rb_3,\rb_4)\in \Lambda_{k_1,k_2,t}^{(u)},\right. \\ & \hspace{4cm} \left.\{\rb_{\pi(1)}\cup \rb_{\pi(2)}\}\cap \{\rb_{\pi(3)}\cup \rb_{\pi(4)}\}=\{t\} \mbox{ for some }\pi\in\mathbb{S}_4 \right\}.  
 \end{aligned}
 \end{equation}
 If $(\rb_1,\rb_2,\rb_3,\rb_4)\in \mathcal{A}$ satisfies \eqref{e:b}, we must have $\rb_{\pi(1)}=\rb_{\pi(2)}$ and $\rb_{\pi(3)}=\rb_{\pi(4)}$ for some $\pi \in \mathbb{S}_4$. Lemma~\ref{l:yprop}(b) then implies
      $$\Ex[Y_{\rb_1}Y_{\rb_2}Y_{\rb_3}Y_{\rb_4}] \lesssim D_{\rb_{\pi(1)}}^2D_{\rb_{\pi(3)}}^2.$$
      where $D_{\rb}$ is defined in~\eqref{d:dw}. However, by definition of $\skr(G_n)$ in~\eqref{d:skgn},
      $$\sum_{t=1}^{v(G_n)}\sum_{(\rb_1,\rb_2,\rb_3,\rb_4)\in \mathcal{A}} \Ex[Y_{\rb_1}Y_{\rb_2}Y_{\rb_3}Y_{\rb_4}]  \le \skr(G_n).$$

 \noindent{\textbf{Case 2.}} We assume either $t_1\neq t_2$ or we have both that $t_1=t_2$ and that $(\rb_1,\rb_2,\rb_3,\rb_4)\not\in \mathcal{A}$ (as defined in~\eqref{adef}). 
 
 \smallskip

Suppose $H_1\cup H_2\cup H_3 \cup H_4$ is a $4$-join of $H$ with $\rb_i\subset V(H_i).$ By~\eqref{e:b}, we observe that $\bigcup_{i=1}^4 H_i$ satisfies~\eqref{good}. We claim that $\bigcup_{i=1}^4 H_i$ is also satisfies \eqref{vgood}. When $t_1=t_2$ and $(\rb_1,\rb_2,\rb_3,\rb_4)\not\in \mathcal{A}$, it is clear from the definition of $\mathcal{A}$, that $\bigcup_{i=1}^4 H_i$ satisfies \eqref{vgood}.

Consider the case when $t_1\neq t_2$. Note that $t_1$ is common to $H_{1},H_{2}$ and $t_2$ is common to $H_{3},H_{4}$. Consequently,  
      $$v\left(\bigcup_{i=1}^4 H_i\right)\leq \min\{v(H_{1}\cup H_{3})+v(H_{2}\cup H_{4})-2, v(H_{1}\cup H_{4})+v(H_{2}\cup H_{3})-2\}.$$
 \noindent\textbf{Case 2.1.} If $u<k_1\wedge k_2$, then there are at least two vertices which appear exactly once among $\rb_1$ and $\rb_2$. By~\eqref{e:b}, these two vertices must be in $\rb_3\cup \rb_4$. Thus, 
      $v(\bigcup_{i=1}^4 H_i)\leq v(H_{1}\cup H_{2})+v(H_{3}\cup H_{4})-2$,
      which forces $\bigcup_{i=1}^4 H_i$ to satisfy \eqref{vgood}, meaning that $\bigcup_{i = 1}^r H_i \in \gj.$

 \medskip 
 \noindent\textbf{Case 2.2.} If $u=k_1\wedge k_2$, then by assumption, $k_1\neq k_2$. Without loss of generality assume $k_1<k_2$. Since $u=k_1$, we have $\rb_1 \subset \rb_2$ and $\rb_3\subset \rb_4$.
 
 The expectation $\Ex[Y_{\rb_1}Y_{\rb_2}Y_{\rb_3}Y_{\rb_4}]$ is nonzero only if $\bigcup_{i = 1}^4 H_i \in \gj$ for any $4$-join of $H$ with $\rb_i \subset V(H_i)$. By way of contradiction, suppose that $v(\bigcup_{i=1}^4 H_i)\ge v(H_{1}\cup H_{2})+v(H_{3}\cup H_{4})-1$. This implies that $|\rb_2 \cap \rb_4|\le 1$.
 If $\rb_2\cap \rb_4=\emptyset$, then $\Ex[Y_{\rb_1}Y_{\rb_2}Y_{\rb_3}Y_{\rb_4}]=\Ex[Y_{\rb_1}Y_{\rb_2}]\Ex[Y_{\rb_3}Y_{\rb_4}]=0$. If $\rb_2\cap \rb_4=\{v\}$ for some vertex $v$, taking the conditional expectation with respect to $X_v$ gives
 $$\Ex[Y_{\rb_1}Y_{\rb_2}Y_{\rb_3}Y_{\rb_4}]=\Ex[\Ex[Y_{\rb_1}Y_{\rb_2}\mid X_v]\Ex[Y_{\rb_3}Y_{\rb_4}\mid X_v]]=\Ex[\Ex[Y_{\rb_1}Y_{\rb_2}]\Ex[Y_{\rb_3}Y_{\rb_4}]].$$
The last equality follows from Lemma~\ref{l:yprop}(c). Since $\Ex[Y_{\rb_1}Y_{\rb_2}]=0$ by Lemma~\ref{l:yprop}(a), the term above is zero. Thus, the expectation is nonzero only if the corresponding $\bigcup_{i=1}^4 H_i$ satisfy~\eqref{vgood}.   

\medskip

\noindent By Lemma~\ref{l:yprop}(b), $\Ex[Y_{\rb_1}Y_{\rb_2}Y_{\rb_3}Y_{\rb_4}] \lesssim D_{\rb_1}D_{\rb_2}D_{\rb_3}D_{\rb_4}$. This expression counts the number of $4$-joins $\bigcup_{i=1}^4 H_i$ of $H$ with $\rb_i\subset V(H_i)$ for $i\in \ll1,4\rr$ (see \eqref{d:dw}). By the argument above, in \textbf{Case 2} it suffices to consider only the $\bigcup_{i=1}^4 H_i$ which lie in $\gj$. Thus, summing over all possibilities in \textbf{Case 2} produces $N(\gj,G_n)$ as an upper bound. Combining this with the estimate from \textbf{Case 1} proves Claim~\ref{e:b1}(b).  
\end{proof}

\section{Proof of Theorem~\ref{t:fourthmoment}}\label{sec4}

{In this section we characterize the fourth moment of $Z(H,G_{n})$ in terms of good joins and we prove Theorem \ref{t:fourthmoment}, the fourth moment result for our central limit theorem.}
For the duration of this section we associate an $r$-clique to $\bs\in V(G_n)_r$, by joining all pairs of vertices $\{s_i,s_j\}$ for $i\neq j$ and we define
$$Z_\bs=\ind\{X_{\bs}\}-c^{1-r}.$$
We begin by characterizing the expectation of products of $Z_{\bs_i}$'s in Lemmas~\ref{l:zexpec} and~\ref{l:not-very-good-factors}. These estimates will be crucial in proving Theorem~\ref{t:fourthmoment} at the end of this section.

\begin{lemma}\label{l:zexpec} 
Let $\s_1,\s_2,\cdots  ,\s_L \in V(G_n)_r$ and set $Z_{\s_i}=\ind\{X_{=\bs_i}\}-{c^{1-r}}$ for $i\in [L]$. Suppose that the $r$-cliques associated with the tuples $\bs_1, \ldots, \bs_L$ collectively form a connected graph.
\begin{enumerate}[(a)]
    \item If for all $i \in [L]$,
    \begin{align}\label{e:gd}
        \left|\bigcup_{j=1}^L \brr{\bs}_j\right|-\left|\bigcup_{j\neq i}^L\brr{\bs}_j\right| \le r-2,
    \end{align}
then for all $c\ge 2^{L+1}-2$ we have 
    \begin{align*}
  \tfrac12c^{1-|\cup_{i=1}^L \brr{\bs}_i|} \le  \EE [Z_{\s_1}Z_{\s_2}\cdots Z_{\s_L}] \le \tfrac32c^{1-|\cup_{i=1}^L \brr{\bs}_i|}.
\end{align*}
\item If $\{\s_1,\ldots,\s_L\}$ does not satisfy~\eqref{e:gd}, then $\EE[Z_{\bs_1} Z_{\bs_2}\ldots Z_{\bs_L}] = 0.$
\end{enumerate}
\end{lemma}
\begin{proof} 
\begin{enumerate}[(a)]
    \item We expand $\EE[Z_{\s_1}Z_{\s_2}\cdots Z_{\s_L}]$ and apply linearity of expectation to find 

\begin{align}\nonumber
    \EE[Z_{\s_1}Z_{\s_2}\cdots Z_{\s_L}] & = \sum_{A \subset [L]} (-1)^{L - |A|} \EE\left[ c^{(1-r)(L-|A|)}\prod_{j\in A} \ind\{X_{=\s_{j}}\}\right] \\ & = \sum_{A \subset [L]} (-1)^{L - |A|} c^{(1-r)(L-|A|)+1-|\bigcup_{j\in A} \brr{\s}_{j}|}. \label{e:p1}
\end{align}

Suppose that $A \neq [L]$. In that case, without loss of generality, suppose that $A = \{1, \ldots, a\}$ for some $0 \le a < L$. Let $F_A = \{\bs_1, \ldots, \bs_a\}$ be the graph formed by joining the $a$ $r$-cliques associated to the tuples $\bs_i$. We consider iterative additions to $F_A$ by adding $\bs_j$ for $j > a$ one at a time until we reach $F_{[L]}$. Since $F_{[L]}$ is connected, we can choose an ordering so that each $\bs_j$ for $j > a$ adds at most $r-1$ unique vertices. Further, due to \eqref{e:gd}, the final tuple $\bs_L$ can add at most $r-2$ unique vertices to  $F_{[L]}$. Therefore,

$$\bigg|\bigcup_{j=1}^L \brr{\s}_{j}\bigg| \le (r-1)(L-|A|)+\bigg|\bigcup_{j\in A} \brr{\s}_{j}\bigg|-1.$$
The $A=[L]$ term in~\eqref{e:p1} evaluates to 	$c^{1-|\bigcup_{j=1}^L \brr{\s}_j|}$. By moving this term to the left hand side of~\eqref{e:p1} and taking absolute values, we find that whenever $c > 2^{L+1} - 2$,
\begin{align*}
\left|\EE[Z_{\bs_1} \cdots Z_{\bs_L}]- c^{1 - |\bigcup_{j=1}^L \brr{\s}_j|}\right| & \le \sum_{A \subsetneq [L]} c^{-|\bigcup_{j=1}^L \brr{\s}_j|} \le (2^L - 1)c^{-|\bigcup_{j=1}^L \brr{\s}_j|} \le \tfrac12 c^{1 - |\bigcup_{j=1}^L \brr{\s}_j|}.
\end{align*}

\item Take $\{\s_1,\s_2,\ldots,\s_L\}$ not satisfying \eqref{e:gd}. Without loss of generality we may assume $|\bigcup_{j=1}^L \brr{\s}_j|=|\bigcup_{j=1}^{L-1} \brr{\s}_j|+r-1.$ That means the $r$-tuple $\s_L$ introduces $r-1$ new vertices, which we will call $u_1,u_2,\ldots,u_{r-1}$. Applying the tower property of expectation and taking conditional expectation over all vertices except $u_1, \ldots, u_{r-1}$, we find
\begin{align*}
\EE[Z_{\s_1}\ldots Z_{\s_L}] &= \EE\left[\Ex[Z_{\s_1}\ldots Z_{\s_L}\mid X_{v}:v\in \cup_{j=1}^L \brr{\s}_j\setminus \{u_1,\ldots,u_{r-1}\}]\right]\\
&= \EE\left[Z_{\s_1}\ldots Z_{\s_{L-1}}\Ex[Z_{\s_L}\mid X_{v}:v\in \cup_{j=1}^L \brr{\s}_j\setminus \{u_1,\ldots,u_{r-1}\}]\right]=0.
\end{align*}
The last expectation is zero as $\Ex[Z_{\s_L}\mid X_v]=0$ for any $v\in \s_L$. This completes the proof.
\end{enumerate}
\end{proof}

\begin{lemma}\label{l:not-very-good-factors}
If $\s_1,\s_2,\s_3,\s_4\in V(G_n)_r$ with $|\brr\s_1\cap \brr\s_2|\ge 2$, $|\brr\s_3\cap \brr\s_4|\ge 2$, and $|\{\brr\s_1\cup \brr\s_2\}\cap \{\brr\s_3\cup \brr\s_4\}|=1$ then
$$\Ex[Z_{\s_1}Z_{\s_2}Z_{\s_3}Z_{\s_4}]=\Ex[Z_{\s_1}Z_{\s_2}]\Ex[Z_{\s_3}Z_{\s_4}].$$
\end{lemma}
\begin{proof}
Suppose $\{v\}=\{\brr\s_1\cup \brr\s_2\}\cap \{\brr\s_3\cup \brr\s_4\}$.  Conditioning on $v$ gives us 
 \begin{align*}
     \Ex[Z_{\s_1}Z_{\s_2}Z_{\s_3}Z_{\s_4}\mid X_v] & = \Ex[Z_{\s_1}Z_{\s_2}\mid X_v]\Ex[Z_{\s_3}Z_{\s_4}\mid X_v] \\ & =(c^{1-|\s_1\cup\s_2|}-c^{2-2r})(c^{1-|\s_3\cup\s_4|}-c^{2-2r})=\Ex[Z_{\s_1}Z_{\s_2}]\Ex[Z_{\s_3}Z_{\s_4}].
 \end{align*}
 We get the desired result by taking the expectation of this expression again.
\end{proof}

We need one more observation in order to conclude Theorem~\ref{t:fourthmoment} from Theorem~\ref{t:clt}.

\begin{proposition}\label{l:Wbound} Let $$W_{\s}:=\frac{1}{|\Aut(H)|}a_{H, \bs}(\ind\{1_{X_{=\s}}\}-{c}^{1-r})=\frac{1}{|\Aut(H)|}a_{H, \bs}Z_{\bs},$$
where $Z_{\bs}:=\ind_{X_{=\s}}-c^{1-r},$ and observe that $Z(H,G_n)=\frac1{\sigma_H}\sum_{\s} W_s.$ For all $c\ge 30$ we have 
	\begin{align}\label{l3}
		\Ex[Z(H,G_n)^4-3] \gtrsim \frac1{\sigma_H^4}\sum_{\{\s_1,\s_2,\s_3,\s_4\} \in \gk} \Ex [W_{\s_1}W_{\s_2}W_{\s_3}W_{\s_4}],
	\end{align}
	where $\gk:=\mathcal{G}_{{K}_r}$ and $\{\bs_1,\bs_2,\bs_3,\bs_4\}\in \gk$ means the sum is over all quadruples of $r$-tuples whose associated four $r$-cliques form a good join of ${K}_r$ (see Definition~\ref{def:good}).
\end{proposition}
\begin{proof} Let $\mathcal{Q}_2^r = \mathcal{Q}_2^{K_r}=\{\{\bs_1,\bs_2\} \mid \bs_i\in V(G_n)_r, |\brr\s_1\cap \brr\s_2|\ge 2\}$. Note that
	\begin{align}\label{e:sig-iden}
	\sigma_H^4 & =(\sigma_H^2)^2  =  	
	\sum_{\{\s_1,\s_2\}\in \mathcal{Q}_2^r}  \sum_{\{\s_1,\s_2\}\in \mathcal{Q}_2^r} \Ex[W_{\s_1}W_{\s_2}]\Ex [W_{\s_3}W_{\s_4}].
\end{align}
We are able to restrict the sum to $|\brr\s_1\cap\brr\s_2|\ge 2$ and $|\brr\s_3\cap\brr\s_4|\ge 2$ due to Lemma~\ref{l:zexpec}(b). On the other hand,
$$\Ex[Z(H,G_n)^4] = \frac1{\sigma_H^4}\sum_{\s_1,\s_2,\s_3,\s_4\in V(G_n)_r}\Ex[W_{\s_1}W_{\s_2}W_{\s_3}W_{s_4}].$$

Given $\{\s_1,\s_2,\s_3,\s_4\}$, $\Ex[W_{\s_1}W_{\s_2}W_{\s_3}W_{s_4}]$ is nonzero in only a few cases: 

\begin{enumerate}
	\item[(a)] $\{\s_1,\s_2,\s_3,\s_4\}$ is connected. By Lemma~\ref{l:zexpec}(b), the expectation is nontrivial only if~\eqref{e:gd} is satisfied. Thus either
	\begin{enumerate}
		\item[(i)] $\{\s_1,\s_2,\s_3,\s_4\}\in \gk$ or
		\item[(ii)]  $\{\s_1,\s_2,\s_3,\s_4\}$  can be written as $F_1 \cup F_2$ where $F_1=\{\s_{\pi_1}\cup \s_{\pi_2}\}$ and $F_2=\{\s_{\pi_3}\cup \s_{\pi_4}\}$ and $|F_1\cap F_2|=1$ for some permutation $\pi \in \mathbb{S}_4$. Note that \eqref{e:gd} forces $F_1, F_2 \in \mathcal{Q}_2^r$. Since there are three ways to form the pairing, the overall contribution of the expectation coming from terms $\Ex[W_{\s_1}W_{\s_2}W_{\s_3}W_{s_4}]$ when $\{\s_1,\s_2,\s_3,\s_4\}$ is connected is 
			$$\sum_{\{\s_1,\s_2,\s_3,\s_4\} \in \gk} \Ex[W_{\s_1}W_{\s_2}W_{\s_3}W_{\s_4}]+3 \underbrace{\sum_{ \{\s_1,\s_2 \} \in \mathcal{Q}_2^r } \sum_{\{\s_3,\s_4 \} \in \mathcal{Q}_2^r }}_{\{\brr \s_1 \bigcup \brr \s_2\} \bigcap \{ \brr \s_3 \bigcup \brr \s_4\}=1} \Ex[W_{\s_1}W_{\s_2}]\Ex[W_{\s_3}W_{\s_4}].$$
	(The expectation in the second term splits up due to Lemma~\ref{l:not-very-good-factors}).
	\end{enumerate}
 
	\item[(b)] $\{\s_1,\s_2,\s_3,\s_4\}$ has two connected components. 
	By Lemma~\ref{l:zexpec}(b), the corresponding expectation is nontrivial only if both the connected components are $2$-shared $2$-joins. Since there are three ways to pair up terms from the quadruple, the overall contribution of terms $\Ex[W_{\s_1}W_{\s_2}W_{\s_3}W_{s_4}]$ when $\{\s_1,\s_2,\s_3,\s_4\}$ has two connected components is
	$$3 \underbrace{\sum_{ \{\s_1,\s_2 \} \in \mathcal{Q}_2^r} \sum_{\{\s_3,\s_4 \} \in \mathcal{Q}_2^r }}_{\{\brr \s_1 \bigcup \brr \s_2\} \bigcap \{ \brr \s_3 \bigcup \brr \s_4\}= \emptyset} \Ex[W_{\s_1}W_{\s_2}]\Ex[W_{\s_3}W_{\s_4}].$$
	Here the expectation factors by independence.
\end{enumerate}
Consequently,
\begin{align}
\notag& \Ex[Z(H,G_n)^4] -3 \\ 
\notag & = \frac1{\sigma_H^4}\sum_{\{\s_1,\s_2,\s_3,\s_4\} \in \gk} \Ex[W_{\s_1}W_{\s_2}W_{\s_3}W_{\s_4}] +\frac3{\sigma_H^4} \underbrace{\sum_{ \{\s_1,\s_2 \} \in \mathcal{Q}_2^r } \sum_{\{\s_3,\s_4 \} \in \mathcal{Q}_2^r }}_{\{\brr \s_1 \bigcup \brr \s_2\} \bigcap \{ \brr \s_3 \bigcup \brr \s_4\}= 1} \Ex[W_{\s_1}W_{\s_2}]\Ex[W_{\s_3}W_{\s_4}] \\ \notag & \hspace{3cm} +\frac3{\sigma_H^4} \underbrace{\sum_{ \{\s_1,\s_2 \} \in \mathcal{Q}_2^r } \sum_{\{\s_3,\s_4 \} \in \mathcal{Q}_2^r }}_{\{\brr \s_1 \bigcup \brr \s_2\} \bigcap \{ \brr \s_3 \bigcup \brr \s_4\}= \varnothing } \Ex[W_{\s_1}W_{\s_2}]\Ex[W_{\s_3}W_{\s_4}] - 3 \\ \notag & = \frac1{\sigma_H^4}\sum_{\{\s_1,\s_2,\s_3,\s_4\} \in \gk} \Ex[W_{\s_1}W_{\s_2}W_{\s_3}W_{\s_4}] -\frac3{\sigma_H^4} \underbrace{\sum_{ \{\s_1,\s_2 \} \in \mathcal{Q}_2^r } \sum_{\{\s_3,\s_4 \} \in \mathcal{Q}_2^r }}_{\{\brr \s_1 \bigcup \brr \s_2\} \bigcap \{ |\brr \s_3 \bigcup \brr \s_4\}|\ge 2} \Ex[W_{\s_1}W_{\s_2}]\Ex[W_{\s_3}W_{\s_4}]. 
\end{align} 
The second equality follows from~\eqref{e:sig-iden}. Whenever $\{\s_1,\s_2 \},\{\s_3,\s_4 \} \in \mathcal{Q}_2^r $ with $\{\brr \s_1 \bigcup \brr \s_2\} \bigcap \{ |\brr \s_3 \bigcup \brr \s_4\}|\ge 2$, the corresponding $4$-join must be a good join of ${K}_r$ (see Definition \ref{def:good}). By Lemma~\ref{l:zexpec} we know that all the individual expectations appearing above are nonnegative for $c\ge 30$. Thus,
	\begin{align}\label{41}
		 \Ex[Z(H,G_n)^4] -3 & \ge  
		  \frac1{\sigma_H^4} \sum_{\{\s_1,\s_2,\s_3,\s_4\} \in  \gk} \bigg\{\Ex[W_{\s_1}W_{\s_2}W_{\s_3}W_{\s_4}]-3\Ex[W_{\s_1}W_{\s_2}]\Ex[W_{\s_3}W_{\s_4}]\bigg\}.
	\end{align} 
Applying Lemma~\ref{l:zexpec} and taking $c\ge 30$, then, if $\{\s_1,\s_2,\s_3,\s_4\} \in \gk$,
$$\Ex[Z_{\s_1}Z_{\s_2}Z_{\s_3}Z_{\s_4}] \ge \tfrac12c^{1-|\s_1\cup \s_2\cup \s_3\cup \s_4|} \ge \tfrac12c \cdot c^{1-|\s_1\cup \s_2|+1-|\cup \s_3\cup \s_4|} \ge \frac{20}{3} \cdot\Ex[Z_{\s_1}Z_{\s_2}]\Ex[Z_{\s_3}Z_{\s_4}].$$
Thus for $c\ge 30$, we have $$\Ex[W_{\s_1}W_{\s_2}W_{\s_3}W_{\s_4}]-3\Ex[W_{\s_1}W_{\s_2}]\Ex[W_{\s_3}W_{\s_4}]\ge \frac12\Ex[W_{\s_1}W_{\s_2}W_{\s_3}W_{\s_4}].$$
Plugging the above bound into~\eqref{41} proves~\eqref{l3}.
\end{proof} 

\begin{proof}[Proof of Theorem \ref{t:fourthmoment}] We use the same notation as in Proposition~\ref{l:Wbound}. Applying Lemma~\ref{l:zexpec} and taking $c\ge 30$, we find that
\begin{align*}
    \sum_{\{\s_1,\s_2,\s_3,\s_4\} \in \gk} \Ex [W_{\s_1}W_{\s_2}W_{\s_3}W_{\s_4}]\gtrsim\sum_{\{\s_1,\s_2,\s_3,\s_4\} \in \gk} a_{H,\bs_1}a_{H,\bs_2}a_{H,\bs_3}a_{H,\bs_4} \Ex [Z_{\s_1}Z_{\s_2}Z_{\s_3}Z_{\s_4}] \gtrsim  N(\gj,G_n).
\end{align*}
To see the final inequality, note that $ a_{H,\bs_i}$ is nonzero only when $G_n$ contains a copy of $H$ on each $\bs_i.$ Consequently, the sum over quadruples that form good joins of $r$-cliques can be reduced to counting good joins of $H$.
We estimate  $\sigma_H^2$ as in Lemma~\ref{sigbd} and then conclude by Proposition~\ref{l:Wbound} that
\begin{align*}
\Ex[Z(H,G_n)^4-3] \gtrsim \frac{N(\gj,G_n)}{N(\mfh,G_n)^{2}}.
\end{align*}
This shows that Theorem~\ref{t:fourthmoment} follows from~\eqref{e:clt}.
\end{proof}

\section{Proof of Theorem~\ref{t:converse}}\label{sec5}
It remains to show Theorem~\ref{t:converse} (that we have convergence of the fourth moment $\mathbb{E}[Z(H, G_n)^{4}] \rightarrow 3$ as soon as asymptotic normality holds for $Z(H, G_n)$). The proof relies on a uniform bound on the higher moments of $Z(H, G_n)$.

\begin{theorem}\label{thm:mom-bd} For all $m\in \Z_{>0}$ and $c\ge 2$ we have that for some $C(c,r,m) < \infty$
	$$\sup_{n\in \Z_{>0}}\Ex\left[\left|Z(H,G_n)\right|^m\right] \le C(c,r,m).$$
\end{theorem}

\begin{proof}[Proof of Theorem~\ref{t:converse} given Theorem~\ref{thm:mom-bd}]
By applying Theorem~\ref{thm:mom-bd} with $m = 8$, we observe that the family $\{Z(H, G_n)^4\}_{n = 1}^{\infty}$ is uniformly integrable, 
\begin{align*}
\lim_{a \rightarrow \infty} \sup_n \EE\left[ Z(H, G_n)^4 \ind_{\{Z(H, G_n)^4 > a\}} \right] &\le \lim_{a \rightarrow \infty} \P(Z(H, G_n)^4 > a) \sup_n \EE\left[Z(H, G_n)^8\right] \\
&\le \lim_{a \rightarrow\infty} C(c,r,8) \P(Z(H, G_n)^4 > a) = 0. 
\end{align*}
Consequently, $\EE[Z(H, G_n)^4] \rightarrow \EE[\mcN(0, 1)^4] = 3$.
\end{proof}

The proof Theorem~\ref{thm:mom-bd} comprises the rest of this section. It will be helpful to recall some basic notions related to set systems and hypergraphs, as in~\cite{friedgut}.

\begin{definition}\label{def:hyper} A \textit{hypergraph} is a set system $F:=(S, \mathcal{B})$, where $S$ is a set and $\mathcal{B}$ is a collection of subsets of $S$. We call the elements of $\mathcal{B}$ \textit{hyperedges}. Given $v\in S$, the \textit{degree} of a vertex in $F$ is 
\begin{align*}
    d_{F}(v):=|\{A \ni v \mid A\in \mathcal{B}\}|, 
\end{align*}
and $d_{\min}(F):=\min_{v\in S} d_{F}(v)$ is the minimum degree of $F$. Hypergraph $F$ is \textit{$s$-uniform} for integer $s \ge 2$ if all elements of $\mathcal{B}$ have cardinality $s$.
\newline \newline 
More generally, an \textit{$s$-uniform multi-hypergraph} is an $s$-uniform hypergraph which may have repeated hyperedges, thus $\mathcal{B}$ is a multiset.

\end{definition}

\begin{ex}
A graph is a $2$-uniform hypergraph and the hyperedges are the edges. In this case, $d_{F}(v)$ is just the usual degree of a vertex $v$.

Let $F=(S,\mathcal{B})$ with $S=\{1, 2, 3, 4\}$ and $\mathcal{B}=\{\{1,3, 4\}, \{2,3,4\}, \{2,3,4\} \}$. $F$ is a $3$-uniform multi-hypergraph, with $d_{F}(1)=1$, $d_{F}(2)=2$, and $d_{F}(3)=d_{F}(4)=3$.
\end{ex}

Now we are ready to prove Theorem~\ref{thm:mom-bd}.

\begin{proof}[Proof of Theorem~\ref{thm:mom-bd}]
For a color $b$ and vertex $v$, define $R_{v}(b) :=\ind_{X_{v}=b}-\frac1c.$ Observe that $\sum_{b=1}^{c} R_{v}(b)=0$. For any $\rb \in \Lambda_{k,*}$, we let $R_{\rb}(b)=\prod_{i=1}^{k} R_{w_{i}}(b)$.

For $\bs \in V(G_n)_r$, let $\mathcal{D}_k(\s)=\{ \rb \in \Lambda_{k,*} \mid \rb \subset \brr\s\}$ be the collection of strictly increasing tuples of $\overline{\s}$ with $k$ elements. By construction, $|\mathcal{D}_k(\s)|=\binom{r}{k}$. Note that
\begin{align*}
	\ind\{X_{=\s}\}-\frac1{c^{r-1}} = \sum_{b=1}^{c} \left(\prod_{k=1}^r \left(\ind_{X_{s_k}=b}-\frac1c+\frac1c\right)-\frac1{c^r}\right)= \sum_{b=1}^c \sum_{k=2}^r \frac1{c^{r-k}}\sum_{\rb \in \mathcal{D}_k(\s)} R_{\rb}(b).
\end{align*}
Without loss of generality, let $m\in 2\Z_{>0}$. From the above observation and the elementary inequality $(\sum_{i=1}^{p} a_i)^m \le  p^{m-1} \sum a_i^m$, we obtain 
\begin{align*}
	(T(H,G_n)-\Ex T(H,G_n))^{m} & = 
	\left[\frac1{|\operatorname{Aut}(H)|}\sum_{\s\in V(G_n)_r} a_{H, \bs}\left[\ind\{X_{=\s}\}-\frac1{c^{r-1}}\right]\right]^m \\ & \le  \sum_{b=1}^c\sum_{k=2}^r \frac{(c(r-1))^{m-1}}{c^{rm-km}}\left[\frac1{|\operatorname{Aut}(H)|}\sum_{\s\in V(G_n)_{r}} a_{H, \bs} \sum_{\rb\in \mathcal{D}_k(\s)} R_{\rb}(b)\right]^m. 
\end{align*}
By fixing $k$ and interchanging the order of the sum, we observe that
\begin{align*}
	\frac1{|\operatorname{Aut}(H)|}\sum_{\s\in V(G_n)_{r}} a_{H, \bs}\sum_{\rb\in \mathcal{D}_k(\s)} R_{\rb}(b) & = \sum_{\rb\in \Lambda_{k,*}} R_{\rb}(b)\frac1{|\operatorname{Aut}(H)|} \sum_{\s: \brr{\s} \supset \rb} a_{H, \bs}  =\sum_{\rb\in \Lambda_{k,*}} R_{\rb}(b)D_{\rb},
\end{align*}
where the last equality follows by definition of $D_{\rb}$ (see~\eqref{d:dw}). 
Therefore, to complete the argument it suffices to show 
\begin{align}\label{e:toshow}
	\Ex \left[\sum_{\rb\in \Lambda_{k,*}} R_{\rb}(b)D_{\rb}\right]^m \lesssim_{m} \operatorname{Var}(T(H,G_n))^{m/2}.
\end{align}
We expand the left hand side of the above equation to get 
\begin{align}\label{e:toshow2}
	\mbox{l.h.s.~of \eqref{e:toshow}} = \sum_{\rb_1,\ldots,\rb_m \in \Lambda_{k,*}} \left[\prod_{i=1}^m D_{\rb_i} \right]\Ex\left[\prod_{i=1}^m R_{\rb_i}(b)\right].
\end{align}
Let $\mathcal{P}_{k,m}$ denote the set of $k$-uniform multi-hypergraphs that have $m$ hyperedges (per Definition~\ref{def:hyper}). 
Fix any $\rb_1,\rb_2,\ldots,\rb_m\in \Lambda_{k,*}$ and let $F\in \mathcal{P}_{k,m}$ be the corresponding $k$-uniform multi-hypergraph. Whenever there exists some $v\in V(F)$ with degree $d_{F}(v)=1$, 
\begin{align*}
	\Ex\left[\prod_{i=1}^m R_{\rb_i}(b)\right]=\Ex\left[R_v(b)\cdot\tfrac1{R_v(b)}\prod_{i=1}^m R_{\rb_i}(b)\right]=\Ex\left[R_v(b)\right]\Ex\left[\tfrac1{R_v(b)}\prod_{i=1}^m R_{\rb_i}(b)\right]= 0,
\end{align*}
by independence. In all other cases the expectation above is trivially bounded by some constant depending on $m$ (and also $r$ and $c$). This implies 
\begin{align}\label{done}
	 \mbox{r.h.s.~of \eqref{e:toshow2}}  & \lesssim_m \sum_{\substack{F\in \mathcal{P}_{k,m} \\ d_{\min}(F)\ge 2}} \sum_{\s\in V(G_n)_{v(F)}}\prod_{(u_1,u_2,\ldots,u_k) \in E(F)} D_{s_{u_1},s_{u_2},\ldots,s_{u_k}},
\end{align}
where we let $D_{\s}:=D_{\brr\s}$ for $D_{\brr\s}$ defined in~\eqref{d:dw}. Applying \rev{Corollary A.1 in \cite{BHA22} with $w(\brr\s):=D_{\brr\s}$ (see also Lemma 3.3 in \cite{friedgut})} we find  
\begin{align*}
\mbox{r.h.s.~of~\eqref{done}} \lesssim_m \sum_{\substack{F\in \mathcal{P}_{k,m} \\ d_{\min}(F)\ge 2}}\left(\sum_{\rb\in\Lambda_{k,*}} D_{\rb}^2\right)^{m/2} \lesssim_m \vr[T(H,G_n)]^{m/2},
\end{align*}
where the last inequality follows from Lemma~\ref{sigbd} and the fact that $|\mathcal{P}_{k,m}|\lesssim_m 1$. This completes the proof of~\eqref{e:toshow}.
\end{proof}

\bibliographystyle{plain}
\bibliography{clt.bib}

\end{document}